\documentclass[12 pt]{article}
\usepackage{amssymb,amsmath,amsfonts,amsthm}
\newcommand\BBP{{\mathbb {P}}}
\newcommand\BBI{{\mathbb {I}}}
\newcommand\BBR{{\mathbb {R}}}
\newcommand\bkR{{\mathbb {R}}}
\newcommand\BBN{{\mathbb {N}}}
\newcommand\BBE{{\mathbb {E}}}
\newcommand\bkE{{\mathbb {E}}}
\newcommand\E{{\mathbb {E}}}

\newtheorem {Theorem}{Theorem}[section]
\newtheorem {Proposition}{Proposition}[section]
\newtheorem {Corollary}{Corollary}[section]

\theoremstyle{definition}

\newtheorem{Remark}{Remark}[section]
\newtheorem {Fact}{Fact}[section]

\newcommand\no{\noindent}
\newcommand\ssk{\smallskip}

\newcommand\beq{\begin{equation}}
\newcommand\eeq{\end{equation}}

\def\no{\noindent}
\def\ssk{\smallskip}

\def\Y{(Y_i)_{i\in{\mathbb N}}}

\begin{document}

\title{Rates of convergence in
the central limit theorem for martingales in the non stationary setting}

\author{J\'er\^ome Dedecker\footnote{J\'er\^ome Dedecker, Universit\'e de Paris, CNRS, MAP5, UMR 8145,
45 rue des  Saints-P\`eres,
F-75006 Paris, France.},
Florence Merlev\`ede \footnote{Florence Merlev\`ede, LAMA,  Univ Gustave Eiffel, Univ Paris Est Cr\'eteil, UMR 8050 CNRS,  \  F-77454 Marne-La-Vall\'ee, France.} and Emmanuel Rio \footnote{Emmanuel Rio, Universit\'e de Versailles, LMV UMR 8100 CNRS,  \ 45 avenue des Etats-Unis, 
F-78035 Versailles, France.}
}

\maketitle

\abstract{In this paper, we give rates of convergence, for minimal distances and for the uniform distance, between the law of partial sums of martingale differences and the
limiting Gaussian distribution.  More precisely, denoting by $P_{X}$ the law of a random variable  $X$ and by $G_{a}$ the normal distribution ${\mathcal N} (0,a)$, we are interested  by giving quantitative estimates for the convergence of $P_{S_n/\sqrt{V_n}}$ to $G_1$, where $S_n$ is the partial sum associated with either martingale differences sequences or more general dependent sequences, and $V_n= {\rm Var}(S_n)$. Applications to linear statistics, non stationary $\rho$-mixing sequences, and sequential dynamical systems are given.}

\medskip


\medskip


\medskip


\section{Introduction and Notations}
\setcounter{equation}{0}

Let $ (\xi_i)_{i\in {\mathbb N}} $ denote a 
 sequence of martingale differences in ${\mathbb L}^2$, with respect to the increasing filtration $ ( {\mathcal F}_i)_{i\in {\mathbb N}} $.  Let $M_n = \sum_{k=1}^n \xi_k$ and $V_n = \sum_{k=1}^n \E (\xi_k^2)$. If 
\beq \label{condTCL} V_n^{-1/2}{\mathbb E} \left (\max_{1 \leq i \leq n} |\xi_i| \right) \rightarrow 0  \ \mbox{ and }  \ V_n^{-1}   \sum_{k=1}^n \xi_k^2 \rightarrow^{{\mathbb P}}  1 \quad \mbox{as $n \rightarrow \infty $,}
\eeq
 then $ V_n^{-1/2} M_n$ converges in  distribution to a standard normal variable (see \cite{GH79}).  Other sets of conditions implying the central limit theorem can be found in \cite{He82}. In particular, under the first part of condition 
 \eqref{condTCL}, its second part is implied by $$V_n^{-1}  \langle M \rangle_n \rightarrow^{{\mathbb P}}  1  \ \text{ as $n \rightarrow \infty $, where } \ \langle M \rangle_n  := \sum_{k=1}^n {\mathbb E} ( \xi_k^2| {\mathcal F}_{k-1} ) \, . $$

  We are interested in bounds on the speed of convergence in this central limit theorem and in particular by giving upper bounds for the
  ${\mathbb L}_1$ and   ${\mathbb L}_{\infty}$ distances defined  respectively as 
  \beq \label{definitionDeltan}
\Delta_{n,1}:=   \Vert F_n -  \Phi  \Vert_1  \ \text{ and } \ \Delta_{n, \infty} := \Vert F_n -  \Phi  \Vert_{\infty}  \, , \eeq
 where $F_n$ is the cdf of $M_n /{\sqrt V_n}$ and ${\Phi}$ is the cdf of a  standard normal variable.  Both of these distances have their own interests. For instance, $ \Delta_{n, \infty}$ provides useful estimates of the quantile $F^{-1}_n (u)$ of $M_n /{\sqrt V_n}$ when $\min (u, 1-u)$ is large enough, whereas the ${\mathbb L}^1$-distance provides estimates of the super quantile (also called the conditional value at risk) as stated in \cite[Theorem 2]{Rio17}.

Concerning the ${\mathbb L}_{\infty}$-distance $ \Delta_{n, \infty}$ for martingales,  several results  have been obtained under different kinds of assumptions. 
 
 One of the first results is due to Heyde and Brown \cite{HB70} and can be stated as follows. For  $p \in ]2, 4 ]$, there exists a positive constant $C_p$ such that for any $n \geq 1$,
\beq \label{BorneHB}
\Delta_{n,\infty}  \leq C_p \Big ( \Vert  V_n^{-1}  \langle M \rangle_n  -1\Vert_{p/2}^{p/2}
+ V_n^{-p/2}  \sum_{k=1}^n {\mathbb E} ( |\xi_k|^p  ) \Big )^{1/(p+1)} \, .
\eeq
This result has been extended to any $p \in (2, \infty )$ by Haeusler \cite{Ha88}. See also Mourrat \cite{Mou}  for an improvement of \eqref{BorneHB} in the bounded case.  If the conditional variances are constant meaning that $ {\mathbb E} ( \xi_k^2| {\mathcal F}_{k-1} ) =  {\mathbb E} ( \xi_k^2)$ a.s. for any $k$, and if
\beq \label{comparep2}
\sup_{i \geq 1} \frac{{\mathbb E} ( |\xi_i|^p  ) } {{\mathbb E} ( |\xi_i|^2 )}  < \infty \, , 
\eeq
 the rates in the central limit theorem in terms of the   ${\mathbb L}_{\infty}$-distance are of order $V_n^{-(p-2)/ (2p+2) }$. For $p=3$ this gives the rate $V_n^{-1/8}$. However in that case, under the additional assumption that there exist two positive constants $\alpha$ and $\beta$ such that for any $i \geq 1$, $  \alpha \leq  {\mathbb E} ( |\xi_i|^2 )   \leq  \beta $, Grams \cite{Grams} proved that the rate is of order $V_n^{-1/4}$ 
(see Theorem 1 in Bolthausen \cite{Bo82}). Even if this rate can appear to be poor compared with the iid case, it cannot be improved without additional assumptions as shown in \cite[Section 6, Example 1]{Bo82}.  More generally, when $p \in (2,3)$, under the same condition on the conditional variances and assuming  \eqref{comparep2}, one can reach the rate $V_n^{-(p-2)/(2p-2)}$ (see our Corollary \ref{Berry-Esseencor}). Again this rate cannot be improved without additional assumptions as shown by our Proposition \ref{CE}.  The paper \cite{EMO} is in this direction. For instance, still in the case where the conditional variances are constant,  Theorem 2  in \cite{EMO} states that  $\Delta_{n, \infty} \leq C V_n^{-1/2} \log n $ provided 
$V_n \leq 4^n$ and  there exists $\gamma >0$ such that $ {\mathbb E} (| \xi_k|^3| {\mathcal F}_{k-1} ) \leq \gamma  {\mathbb E} ( \xi_k^2| {\mathcal F}_{k-1} )$ a.s. for any $k$ (see \cite{Fan19} for related results).

Let us now comment on the quantity $\Vert  V_n^{-1}  \langle M \rangle_n  -1\Vert_{p/2}$ appearing in the right hand side of \eqref{definitionDeltan} when it is not equal to zero. For stationary sequences (except in some degenerate cases),   $\Vert  V_n^{-1}  \langle M \rangle_n  -1\Vert_{p/2} $ is typically of order $V_n^{-1/2} $ which leads at best  to the rate  $V_n^{-p/(4p+4)}$.  It is therefore clear that, in these non-degenerate situations, the rate  $V_n^{-1/4}$ cannot be reached, whatever the value of $p$.

One of the goals of this paper is to give tractable conditions   (not assuming that ${\mathbb E} ( \xi_k^2| {\mathcal F}_{k-1} ) =  {\mathbb E} ( \xi_k^2)$ a.s. or  $ V_n^{-1}  \langle M \rangle_n  =1 $ a.s.)  for $p \in (2, 3 ] $ under which the rate $V_n^{-(p-2)/(2p-2)}$ can be reached 
for $\Delta_{n, \infty}$ (up to a logarithmic term when $p=3$).  These conditions will be expressed with the help of 
quantities involving a sum of conditional expectations and  allow to use martingale approximations techniques, as introduced by Gordin \cite{Go69} (see  also Voln\'y \cite{Vo93}), to get rates when the sequence is not a martingale differences sequence.  Applications via martingale approximations are provided  in Section \ref{sectionapplications}. The case of sequential dynamical systems as developed by Conze and Raugi \cite{CR} is considered in Subsection \ref{subsectionseqdyn}. 

To derive the rates concerning  $\Delta_{n,\infty}$, we shall rather work with  minimal distances also called Wasserstein distances of order $r$  (see Inequality \eqref{ineBEW1} below for the connection between  $\Delta_{n,\infty}$ and these distances). In particular, we shall  also exhibit rates for the minimal distance $\Delta_{n,1}$ (see the equality \eqref{linkW1Delta1} below). 

Let us  recall the definitions of these minimal distances. Let ${\mathcal L}(\mu, \nu)$ be the set of  probability laws on
$\mathbb R^2$ with  marginals $\mu$ and $\nu$. Let us consider the
following minimal distances: for any $r >0$, 
\[
W_r(\mu, \nu) = \displaystyle \inf \Big \{ \Big (\int |x-y|^r P(dx, dy) \Big )^{1/\max(1,r)}
: P \in {\mathcal L}(\mu, \nu) \Big \}  \, .  \]
%
We consider also the following ideal distances of order $r$
(Zolotarev distances of order $r$). For two probability measures
$\mu$ and $\nu$, and $r$ a positive real, let
$$
\zeta_{r} ( \mu , \nu )   =  \sup \Big \{ \int fd\mu - \int f d \nu
: f \in \Lambda_r \Big \} \, ,
$$
where $\Lambda_r$ is defined as follows: denoting by $l$ the natural
integer such that $ l < r \leq l +1$, $\Lambda_r$ is the class of
real functions $f$ which are $l$-times continuously differentiable
and such that \beq \label{defgamrrio} |f^{(l)}(x) - f^{(l)}(y) |
\leq | x - y |^{r-l} \ \hbox{ for any } (x,y) \in \bkR \times \bkR
\, . \eeq
For $r \in ]0,1]$, applying the  Kantorovich-Rubinstein theorem (see for instance \cite[Theorem 11.8.2]{Du})
to the metric $d(x,y)=|x-y|^r$, we infer that
\beq \label{KRdual} W_r(\mu, \nu) =
\zeta_{r} ( \mu , \nu ) \, . \eeq
For $r>1$ and for probability laws on the real line, the following inequality holds \beq \label{riolink} W_r(\mu, \nu) \leq c_r \big ( \zeta_{r} ( \mu ,
\nu ) \big )^{1/r} \, , \eeq
where $c_r$ is a constant depending
only on $r$ (see \cite[Theorem 3.1]{Rio09}).  Note that for $r=1$, \eqref{KRdual}  ensures that 
\beq \label{linkW1Delta1}
W_1 ( P_{M_n/{\sqrt V_n}} , G_1)  = \zeta_1 ( P_{M_n/{\sqrt V_n}} , G_1) = \Delta_{n,1} \, , 
\eeq
where $P_{M_n/{\sqrt V_n}}$ is the law of $M_n/{\sqrt V_n}$ and $G_{1}$ the  ${\mathcal N} (0,1)$ distribution.

The paper is organized as follows. In Section \ref{sectionmainresults}, we give rates in terms of Zolotarev and then in terms of Wasserstein distances between the law of the martingale having a moment of order $p \in (2,3]$ and  the Gaussian distribution with the same variance. Upper and lower bounds for the uniform distance $\Delta_{n, \infty}$ are provided in Section \ref{sectionBEplusCE}. Applications to linear statistics associated with stationary sequences, $\rho$-mixing sequences in the sense of Kolmogorov and Rozanov \cite{KR} and sequential dynamical systems are presented in Section \ref{sectionapplications}. All the proofs are postponed to Section \ref{sectionproofs}. 

\smallskip

In the rest of the paper, we shall use the following notations:  we will denote by $P_{X}$ the law of a r.v.  $X$ and by $G_{a}$ the  ${\mathcal N} (0,a)$ distribution, and for two sequences $(a_n)_{n \geq 1}$ and $(b_n)_{n \geq 1}$ of positive reals, $a_n \ll b_n$ means there exists a positive constant $C$ not depending on $n$ such that $a_n \leq C b_n$ for any $n\geq 1$. Moreover, given a filtration ${\mathcal F}_{\ell}$, we shall often use the notation $\E_{\ell} ( \cdot) = \E ( \cdot |{\mathcal F}_{\ell}) $.

\section{Rates for Zolotarev and Wasserstein distances} \label{sectionmainresults}

\setcounter{equation}{0}

In this section $ (\xi_i)_{i\in {\mathbb N}} $ will denote a 
 sequence of martingale differences in ${\mathbb L}^2$, with respect to the increasing filtration $ ( {\mathcal F}_i)_{i\in {\mathbb N}} $ and with  ${\mathbb E}(\xi_i^2)=\sigma_i^2$.  We shall use the following notations: 
 \[M_n= \sum_{i=1}^n \xi_i \text{ , }V_n = \sum_{i=1}^n \sigma_i^2 \text{ , } \delta_n =  \max_{ 1 \leq i \leq n} |\sigma_{i}|  \text{ , } v_n(a) = a^2 \delta_n^2 + \alpha V_n  \, , \]
where  $a$ is a positive real and $\alpha = (1+a^2)/a^2$. Moreover, for $p \geq 2$ and $\ell \geq 2$,  we denote by
 \begin{equation} \label{defU}
 U_{\ell,n} (p)  = \Big \Vert( |  \xi_{\ell -1} | \vee \sigma_{\ell-1} )^{p-2} \Big | \sum_{k = \ell}^n ( \E_{\ell-1} ( \xi_k^2) - \sigma_k^2) \Big |   \Big \Vert_1
 \, .
\end{equation}

\begin{Theorem} \label{Th1first}   Let $p \in ]2,3]$ and $r \in (0, p]$. There exist  positive constants $C_{r,p}$ depending on $(r,p)$ and $\kappa_r$ depending on $r$ such that for every positive integer $n$ and any $a \geq 1$,
\begin{multline} \label{boundzetarth1}
\zeta_r(P_{M_{n}}, G_{V_n})  \leq C_{r,p} \Bigl (     \delta_n^r    \int_{a}^{\sqrt{ v_n (a)  / \delta_n^2}} \frac{1}{x^{3-r}} dx   + \delta_n^{r-1}   \int_{a}^{\sqrt{ v_n (a) / \delta_n^2}}  \frac{\psi_n ( \kappa_r  x ) }{x^{2-r}}dx  +L_n(p,r,a \delta_n)   \Bigr )  \\ + 4\sqrt{2}  a^r \delta_n^r
 \, ,
\end{multline}
where 
\begin{equation} \label{defpsi}
 \psi_n(t) = \sup_{1 \leq k \leq n } \frac{ \E \inf ( t  \delta_n \xi_k^2 , |\xi_k|^3)}{\sigma_k^2}  
\end{equation}
and 
\begin{equation} \label{defLn}
L_n(p,r,a \delta_n) = \sum_{ \ell =2}^n \frac{U_{\ell,n} (p) }{ (V_n - V_{\ell -1}  + a^2  \delta_n^2)^{(p-r)/2}}   \, .
\end{equation}
\end{Theorem}

\begin{Remark} \label{remBErate}
Let $p \in ]2,3]$ and $r \in (0, p]$. Using \eqref{KRdual} or  \eqref{riolink}, the fact that   
$$
\zeta_r(P_{M_{n}/\sqrt{V_n}}, G_{1})  = V_n^{-r/2} \zeta_r(P_{M_{n}}, G_{V_n}) 
$$
 and inequality \eqref{boundzetarth1}, we derive upper bounds for 
$W_r(P_{M_{n}/\sqrt{V_n}}, G_{1}) $ and then rates in the central limit theorem. In particular for $W_r(P_{M_{n}/\sqrt{V_n}}, G_{1}) $ to converge to zero as $n \rightarrow \infty$ it is necessary that 
$V_n^{-1/2}  \max_{ 1 \leq i \leq n} |\sigma_{i}| \rightarrow 0$ as $n \rightarrow \infty$ which is also a necessary  condition for the CLT to hold. 
\end{Remark}

In particular, for $r \in (0,1]$, the following corollary holds.
\begin{Corollary} \label{corW1}
 Let $p \in ]2,3]$ and $r \in (0,1]$. Under the assumptions and notations of Theorem \ref{Th1first}, there exists a positive constant $C_{r,p}$ depending on $(r,p)$ such that 
\begin{equation*} \label{boundW1}
W_r(P_{M_{n}}, G_{V_n})  \leq 4\sqrt{2}  (a \delta_n)^r + C_{r,p} \left (        \int_{a}^{\sqrt{ v_n(a)  / \delta_n^2}}  \frac{\psi_n (6x ) }{x}dx +L_n(p,r,a \delta_n)  \right )  \, .
\end{equation*}
In particular if the $\xi_i$'s are 
 in ${\mathbb L}^p$ with $p \in ]2,3]$ and $(r,p) \neq (1,3)$, 
  \begin{equation*} \label{boundW1-pdiff3}
W_r(P_{M_{n}}, G_{V_n})  \leq 4\sqrt{2}  (a \delta_n)^r + { \tilde C}_{r,p} \left (     \sup_{1 \leq k \leq n } \frac{ \E( |\xi_k|^p)}{\sigma_k^2} ( v_n (a) )^{(2+r-p)/2} + L_n(p,r,a \delta_n)    \right )
 \, ,
\end{equation*}
 and if the $\xi_i$'s are  
 in ${\mathbb L}^3$, 
 \begin{equation*} \label{boundW1-p=3}
W_1(P_{M_{n}}, G_{V_n})  \leq 4\sqrt{2}  a \delta_n + { \tilde C}_{3} \left (     \sup_{1 \leq k \leq n } \frac{ \E( |\xi_k|^3)}{\sigma_k^2} \log ( \sqrt{ v_n (a) } / \delta_n) + L_n(3,1,a \delta_n)     \right )
 \, .
\end{equation*}
\end{Corollary}
\begin{Remark} Note that if $(\xi_i)_{i \geq 1}$ is a sequence of  integer valued random variables   then, whatever its dependence structure, setting $S_n = \sum_{k=1}^n \xi_i$ and proceeding as in the proof of \cite[Theorem 5.1]{Rio09} we derive that for any $r >0$, 
\[
\liminf_{n \rightarrow \infty}   \Big (  W_r(P_{S_{n}}, G_{{\rm Var} (S_n)})    \Big )^{\max(1,r)} \geq 2^{-r}/(r+1) 
\]
provided ${\rm Var} (S_n) \rightarrow  \infty$ as $n \rightarrow \infty$. Hence, in the case of martingale differences, if $p \in (2,3)$, $\sup_{1 \leq k \leq n } \sigma_k^{-2} \E( |\xi_k|^p) \leq C_1$ and  $L_n(p,p-2,  \delta_n) \leq C_2$, we get 
$$
  2^{-(p-2)}/(p-1) \leq   \liminf_{n \rightarrow \infty} W_{p-2}(P_{M_{n}}, G_{V_n}) \leq \limsup  _{n \rightarrow \infty} W_{p-2}(P_{M_{n}}, G_{V_n}) \leq K
$$
for some positive constant $K$. In addition, if $p=3$, $\sup_{1 \leq k \leq n } \sigma_k^{-2}\E( |\xi_k|^3) \leq C_1$ and  $L_n(3,1,  \delta_n) \leq C_2$, we have 
$$
  W_{1}(P_{M_{n}}, G_{V_n})  \ll \log ( \sqrt{ v_n (1) } / \delta_n)  \, .
$$

\end{Remark}

\section{Berry-Esseen type results}  \label{sectionBEplusCE}
\setcounter{equation}{0} 
Using  \cite[Remark 2.4]{DMR} stating that, for any $p \in ]2,3]$ and any integrable real-valued random variable $Z$, 
\beq \label{ineBEW1}
\sup_{x \in {\mathbb R}} \big |  {\mathbb P} ( Z \leq x ) - \Phi (x)  \big |  \leq  ( 1+ (2 \pi)^{-1/2} ) \big (W_{p-2}(P_{Z}, G_{1})  \big )^{1/(p-1)} \, ,
\eeq
combined with Remark \ref{remBErate},  Corollary \ref{corW1} leads also to Berry-Esseen type upper bounds.  More precisely, the following result holds
\begin{Corollary} \label{Berry-Esseencor} Assume that $(\xi_i)_{i\in {\mathbb Z}}$ is a sequence of martingale differences 
 in ${\mathbb L}^p$ with $p \in ]2,3]$. Let $\Delta_{n, \infty}$ be defined by \eqref{definitionDeltan}. Then, with the notations of Section \ref{sectionmainresults},  one has 
 \[ \Delta_{n,\infty} \ll      \left\{
\begin{aligned}
& V_n^{- \frac{(p-2)}{2(p-1)}}  \left (     \sup_{1 \leq k \leq n } \frac{ \E( |\xi_k|^p)}{\sigma_k^2} + L_n(p,p-2, \delta_n)    \right )^{1/(p-1)}   & \text{ if $p\in (2,3) $}  \\  
&V_n^{-1/4} \left (     \sup_{1 \leq k \leq n } \frac{ \E( |\xi_k|^3)}{\sigma_k^2} \log (  \sqrt{ v_n (1) } / \delta_n) + L_n(3,1, \delta_n)     \right )^{1/2}   & \text{ if $p=3 $.} \\
  \end{aligned}
\right.
\]
\end{Corollary}
In particular if  
\beq \label{conditionforte}
  \sup_{1 \leq k \leq n } \frac{ \E( |\xi_k|^p)}{\sigma_k^2} \leq C  \quad \mbox{ and } \quad  \E (\xi_k^2 | {\mathcal F}_{k-1} ) = \sigma_k^2 \mbox{ a.s. }
\eeq
it follows that 
\[
\Delta_{n,\infty}   \ll  \left\{
\begin{aligned}
&  V_n^{- \frac{(p-2)}{2(p-1)}}  & \text{ if $p\in (2,3) $}  \\  
&V_n^{-1/4}  \log^{1/2} (  \sqrt{ v_n (1) } / \delta_n) & \text{ if $p=3 $.} 
  \end{aligned}
\right.
\]
It turns out that one can construct a non stationary sequence of martingale differences satisfying \eqref{conditionforte} with $\sigma_k^2=1$ and such that there exists a positive constant $c >0$ for which $\Delta_n \geq c n^{- \frac{(p-2)}{2(p-1)}} $   for any $p >2$ and  any $n\geq 20$. This shows that for $p \in (2,3)$ the rate given in Corollary \ref{Berry-Esseencor} is optimal and quasi optimal (up to $\sqrt{\log n}$)  in case $p=3$. 

\begin{Proposition} \label{CE} Let $p >2$ and $n \geq 20$. There exists $(X_1, \ldots, X_n)$ such that 
\begin{enumerate}
\item $\E(X_k | \sigma( X_1, \dots , X_{k-1} ) )=0 $ and $\E(X^2_k | \sigma( X_1, \dots , X_{k-1} )) =1  $ a.s.,
\item $ \sup_{1 \leq k \leq n} \E(|X_k|^p ) \leq \E ( |Y|^p) + 5^{p-2} $ where $Y \sim {\mathcal N} (0,1) $, 
\item $ \sup_{t\in \BBR}  \big| \BBP ( S_n \leq t \sqrt{n} ) - \Phi (t) \big| \geq 0.06 \,\, n^{-(p-2)/(2p-2)}$, where $S_n= \sum_{k=1}^n X_k$. 
\end{enumerate}
\end{Proposition}

Note that in case $p=3$, Example 1 in \cite{Bo82} also shows that even for martingales with conditional variances equal to one and moments of order  $3$ uniformly bounded, the rate $n^{-1/4}$ cannot be improved in general.  

\medskip

\noindent {\it Proof of Proposition \ref{CE}.}  Let $n$ be  an integer satisfying $n\geq 20$. Let $a$ be a real in $[1, \sqrt{n} /4 [$, to be fixed later, and $k = \inf \{ j \in \BBN : j \geq 4a^2 \}$.
Then $k < 1 + (n/4)$, which ensures that $k<n$.  Set $m=n-k$.  We now define the sequence $(X_j)_{j \in [1,n]}$ of martingale differences
as follows. 
\par\ssk\no
{\bf (i)}  The random variables $(X_j)_{j\in [1,m]}$ are independent and identically distributed with common law the standard normal law.
\par\ssk\no
{\bf (ii)} Let $U_{m+1}, \ldots , U_n$ be a sequence of independent random variables with uniform distribution over $[0,1]$, independent of 
$(X_1, X_2, \ldots , X_m)$.  Let $S_m = X_1 + X_2 + \cdots + X_m$. If $|S_m| \notin [a,2a]$, set $X_j = \Phi^{-1} (U_j)$ for any 
$j$ in $[m+1,n]$. If $|S_m| \in [a, 2a]$, set 
\beq \label{DefXj}
X_j = -(S_m/k) \BBI_{U_j \leq k^2/(S_m^2 + k^2) }  +  (k/S_m) \BBI_{U_j > k^2/(S_m^2 + k^2) } . 
\eeq
\par\ssk 
From the definition of the random variables $X_j$,  if $|S_m| \in [a,2a]$ and $U_j \leq k^2 / ( S_m^2 + k^2)$ for any $j$ in $[m+1,n]$, 
then $S_n=0$. It follows that 
\beq \label{lowerboundAtom}
\BBP ( S_n=0)  \geq \exp \bigl( - k \log ( 1+ 4a^2/k^2) \bigr) \frac{2}{\sqrt{2\pi m}} \int_a^{2a}  \exp ( - x^2/2m) dx . 
\eeq
\par
We now estimate the conditional moments of the random variables $X_j$ for $j>m$.  From the definition of these random variables, 
for any measurable function $f$ such that $f(X_j)$ is integrable
\beq 
\label{CondExpectationf(Xj)}
\BBE ( f (X_j) \mid {\cal F}_{j-1} ) = \BBE ( f (X_j) \mid S_m ) . 
\eeq
Now, if $(S_m = x)$ for some $x$ such that $|x| \notin [a,2a]$, then $X_j = \Phi^{-1} (U_j)$ and consequently 
\beq \label{CondMoments1}
\BBE (X_j \mid S_m=x) = 0 \ , \  \BBE ( X_j^2 \mid S_m = x) = 1 \ \text{ and }  \BBE ( |X_j|^p \mid S_m = x ) = \BBE ( |Y|^p )  
\eeq
for any $p>0$. Here $Y$ is a random variable with law ${\mathcal N}(0,1)$. Next, if $(S_m = x)$ for some $x$ such that $|x| \in [a,2a]$, then,
according to (\ref{DefXj}), 
\beq \label{CondMoments2}
\BBE (X_j \mid S_m=x) = 0 \ , \  \BBE ( X_j^2 \mid S_m = x) = 1 
\eeq
and, for any $p>2$, 
\beq \label{CondMoments3}
\BBE ( |X_j|^p \mid S_m = x ) =  \frac{ |x|^p k^{2-p} + k^p |x|^{2-p}  }{x^2 + k^2}  . 
\eeq
In that case, since $k \in [4a^2 ,5a^2]$ and $|S_n| \in [a,2a]$, 
\beq \label{UpperBoundMoments1}
\BBE ( |X_j|^p \mid S_m = x ) \leq |x|^p k^{-p} + k^{p-2} |x|^{2-p} \leq 1+ (5a)^{p-2} \leq 2 \, (5a)^{p-2}  \, .
\eeq
From (\ref{CondMoments1}), the above upper bound and the fact that, since $n \geq 20$,  $m\geq (3n/4) - 1 \geq (7n/10)$ and then 
\beq \label{UpperBoundMoments2}
\BBE ( |X_j|^p  ) \leq \BBE ( |Y|^p )  + 2\, (5a)^{p-2} \BBP ( |S_m| \in [a,2a] ) \leq  \BBE ( |Y|^p )  +5^{p-2} 2 a^{p-1} n^{-1/2}  \, .
\eeq
Now, for $p>2$, choosing $a = (n/4)^{1/(2p-2)}$ in the above inequality, we get that 
\beq \label{UpperBoundMoments3}
\BBE ( |X_j|^p  ) \leq \ \BBE ( |Y|^p )  +5^{p-2}  \, .
\eeq
Consequently, for this choice of $a$, the absolute moments of order $p$ of the random variables $X_j$ are bounded by some positive constant
 depending only on $p$.  
\par\ssk
Now, using (\ref{lowerboundAtom}) we bound from below $\BBP (S_n=0)$. First $4a^2 \leq k$, which ensures that 
$\exp \bigl( - k \log ( 1+ 4a^2/k^2) \bigr) \geq 1/e$, and second, for $x$ in $[a,2a]$, 
$$
\exp ( -x^2 /2m) \geq \exp ( - 2a^2 / m ) \geq \exp ( - n/8m) \geq \exp (-10/56)  
$$
since $a^2 \leq n/16$ and $m\geq 7n/10$. Hence
\beq \label{lowerboundAtom2}
\BBP ( S_n = 0 ) \geq   0.24 \,\,  a n^{-1/2} \geq 0.12 \,\, n^{-(p-2)/(2p-2)}  \, .
\eeq
Therefrom, Item 3 of the proposition follows.  $\square$

\section{Applications} \label{sectionapplications}

Proposition \ref{Mainprop} of Section  \ref{sectionproofs} (which is the main ingredient for proving Theorem \ref{Th1first}), combined with a suitable martingale approximation, can also be used to derive upper bounds for the Wasserstein distances between the law of  partial sums of non necessarily stationary sequences and the corresponding limiting Gaussian distribution.  This leads to new results for linear statistics, $\rho$-mixing sequences and  sequential dynamical systems. Note that for these non stationary dynamical systems, a reversed martingale version of our Theorem \ref{Th1first}  will be needed.

\setcounter{equation}{0}

\subsection{Linear statistics} Let $p \in ]2,3]$ and $(Y_i)_{i \in {\mathbb Z}}$  be a strictly stationary sequence of centered real-valued random variables in ${\mathbb L}^p$. Let ${\mathcal G}_k = \sigma (Y_i, i \leq k)$. Define  
$\gamma_k = {\rm Cov }(Y_0,Y_k)$ and 
\[
\lambda_k =  \max \Big ( \Vert Y_0 \E(Y_k | {\mathcal G}_0) \Vert_{p/2} ,  \sup_{j\geq i \geq k} \Vert  \E (Y_i Y_j  | {\mathcal G}_0)) - \E ( Y_i Y_j)  \Vert_{p/2}   \Big ) \, .
\]
Let also 
\beq \label{defLambdaeta}
\Lambda_n = \sum_{i =1}^n i \lambda_i  \ \mbox { and }  \  \eta_n =  \sum_{i =0}^n  \Vert \E(Y_i | {\mathcal G}_0) \Vert_p \, .
\eeq

Let $(\alpha_{i,n})_{ i \geq 1} $ a triangular array of real numbers and define 
\[m_n = \max_{1 \leq \ell \leq n} | \alpha_{\ell,n} | \, , \, X_{i,n} = \alpha_{i,n} Y_i   \, , \, S_n = \sum_{i=1}^n X_{i,n} \quad \text{and} \quad V_n = {\rm Var} (S_n) \, . \]

We refer to $S_n$ as a ``linear statistic" based on the stationary sequence $(Y_i)_{i \in {\mathbb Z}}$. Such linear statistics appear in many statistical contexts, for instance when considering least square estimators in a regression model with stationary errors (see for instance \cite{DM}). 

\smallskip

In the two corollaries below we shall assume that $\sum_{k \geq 0} |\gamma_k| < \infty$ which implies in particular that $(Y_i)_{i \in {\mathbb Z}}$ has a bounded spectral density    
$f_Y ( \theta ) = \frac{1}{2 \pi} \sum_{k \in {\mathbb Z}} \gamma_k 
{\mathrm e}^{{\rm i} k \theta} $  on  $[-\pi,\pi]$. Moreover, in the first corollary, we assume in addition  that the spectral density is bounded away from $0$ (we refer to \cite{BR02} for conditions ensuring such a fact).  To state these corollaries, it is  convenient to introduce the 
following quantity: 
\beq \label{defBnp}
B(n,p) :=  \left\{
\begin{aligned}
m_n^{p-2} \eta_n^{p-2} (  \Lambda_n  + \eta_n^2)   \Big (   \sum_{\ell= 1}^{n }\alpha^2_{\ell,n}  \Big  )^{(3-p)/2}  & \  \text{ if $p\in (2,3) $}  \\  
m_n  \eta_n (  \Lambda_n  + \eta_n^2)  \log \Big (  m_n^{-1} \sum_{\ell= 1}^n \alpha^2_{\ell,n}  \Big  )   & \  \text{ if $p=3 $.} \\
  \end{aligned}
\right.
\eeq

\begin{Corollary} \label{corapplilinearcondonspectral}  Let $p\in (2,3] $. Assume that $\sum_{k \geq 0} |\gamma_k| < \infty$ and that   $\inf_{t \in [-\pi, \pi ]} |f_Y(t)| =m >0$. Then 
\[ W_1(P_{S_{n}}, G_{V_n})  \ll   m_n  \sum_{k=0}^n \Vert \E (Y_k | {\mathcal G}_0) \Vert_2   +  B(n,p)  \, .\]
\end{Corollary}
Note that if \beq \label{corL2}
\sum_{i \geq 1} \Vert \E(Y_i | {\mathcal G}_0) \Vert_2 < \infty \, ,
\eeq
then  $\sum_{k \geq 0} |\gamma_k| < \infty$ (see for instance \cite[p. 106]{MPU}). If in addition to \eqref{corL2}, we assume that  $\sup_{n \geq 0}  ( \Lambda_n  + \eta_n) < \infty$, then we get 
\begin{equation} \label{linstat2}
W_1(P_{S_{n}}, G_{V_n})  \ll  \left\{
\begin{aligned}
m_n^{p-2}  \Big (   \sum_{\ell= 1}^{n }\alpha^2_{\ell,n}  \Big  )^{(3-p)/2}  & \ \text{ if $p\in (2,3) $}  \\  
m_n   \log \Big (  m_n^{-1} \sum_{\ell= 1}^n \alpha^2_{\ell,n}  \Big  )   & \  \text{ if $p=3 $.} \\
  \end{aligned}
\right.
\end{equation}
For additional results in the  special case where $(Y_i)_{i \in {\mathbb Z}}$ is a stationary sequence of martingale differences, we refer to \cite{DM}.
\begin{Remark} \label{explicit} 
If, for any positive $k$, 
\[
\lim_{n \rightarrow \infty}  \frac{ \sum_{\ell= 1}^{n-k} \alpha_{\ell,n}  \alpha_{\ell +k,n} }{ \sum_{\ell= 1}^n \alpha^2_{\ell,n} } = c_k \, , 
\]
and $\sum_{k \geq 0} |\gamma_k| < \infty$, then 
\beq \label{condHannanvar}
\frac{V_n}{ \sum_{\ell= 1}^n \alpha^2_{\ell,n} } \rightarrow \sigma^2 = \gamma_0 + 2 \sum_{k \geq 1 } c_k \gamma_k \, , \,  \mbox{ as } n \rightarrow \infty \, .
\eeq
Moreover if $\inf_{t \in [-\pi, \pi ]} |f_Y(t)| =m >0$, then $\sigma^2 >0$.  Let $T_n = S_{n}/ \sqrt{\sum_{\ell= 1}^n \alpha^2_{\ell,n} }$. Under \eqref{corL2} and if $f_Y$ is bounded away from zero, $\sup_{n \geq 0}  ( \Lambda_n  + \eta_n) < \infty$ and \eqref{condHannanvar} holds, it follows that 
\[ W_1(P_{T_{n}}, G_{\sigma^2})  \ll    \left| \frac{ V^{1/2}_n}{ \sqrt{ \sum_{\ell= 1}^n \alpha^2_{\ell,n} }}- \sigma  \right | +  \left\{
\begin{aligned}
&  \left (  \frac{m_n  }{\sqrt{\sum_{\ell= 1}^{n }\alpha^2_{\ell,n} }} \right )^{p-2}   & \text{ if $p\in (2,3) $}  \\  
& \frac{m_n  }{\sqrt{\sum_{\ell= 1}^{n }\alpha^2_{\ell,n} }} \log \Big ( m_n^{-1} \sum_{\ell= 1}^n \alpha^2_{\ell,n}  \Big  )   & \text{ if $p=3 $.} \\
  \end{aligned}
\right.
\]
In case  where $\alpha_{k,n} = \kappa k^{\alpha}$ with $\alpha >-1/2$,  then $m_n ( \sum_{\ell= 1}^{n }\alpha^2_{\ell,n})^{-1/2}$ is exactly of order $n^{-(\alpha +1/2)} {\bf 1}_ {-1/2 < \alpha < 0} + n^{-1/2} {\bf 1}_ {\alpha \geq  0} $ and we can show (since 
$\sum_{i \geq 1} i |\gamma_i| < \infty$ and $\sigma>0$),  that 
\[
 \left| \frac{ V^{1/2}_n}{ \sqrt{ \sum_{\ell= 1}^n \alpha^2_{\ell,n} }}- \sigma  \right | = O (1/n) \, .
\]
Hence, for instance if $\alpha  \geq 0$, 
\[ W_1(P_{T_{n}}, G_{\sigma^2})  \ll     \left\{
\begin{aligned}
&  n^{-(p-2)/2}  & \text{ if $p\in (2,3) $}  \\  
& n^{-1/2} \log (n)  & \text{ if $p=3 $.} \\
  \end{aligned}
\right.
\]
\end{Remark}

\begin{Remark} \label{alphamixing} 
Let $(\alpha_{{\bf Y}} (k) )_{k >0}$ be the usual Rosenblatt strong mixing coefficients \cite{Ros56}  of the sequence $(Y_i)_{i \in {\mathbb Z}}$. 
If we assume that 
\[
 { \mathbb P} ( |Y_0| \geq t ) \leq C t^{-s} \mbox{ for some $s>p $ and  } \sum_{k \geq 1} k ( \alpha_{{\bf Y}} (k) )^{2/p - 2/s} < \infty \, , 
\]
then condition \eqref{corL2} holds and $\sup_{n \geq 0}  ( \Lambda_n  + \eta_n) < \infty$. Hence in this case \eqref{linstat2} holds and  Remark \ref{explicit} applies. 
\end{Remark}

If we do not require the spectral density bounded away from $0$ but only that $f_Y(0) >0$ then an additional term appears in the bound of the Wasserstein distance between $P_{S_{n}}$ and $G_{V_n}$. 
\begin{Corollary} \label{corapplilinear}  
Let $p\in (2,3] $. Assume that  $\sum_{k \geq 1}k^2 | \gamma_k |  < \infty$ and $f_Y(0) >0$. Then 
\[ W_1(P_{S_{n}}, G_{V_n})  \ll  m_n  \sum_{k=0}^n \Vert \E (Y_k | {\mathcal G}_0) \Vert_2 + B(n,p) + \Big (   \sum_{k=1}^{n+1}  ( \alpha_{k,n} -  \alpha_{k-1,n} )^2  \Big )^{1/2} \, , \]
where $B(n,p) $ is defined in \eqref{defBnp}. 
\end{Corollary}

\subsection{$\rho$-mixing sequences} In this section we consider a sequence $(X_i)_{i  \geq 1}$ of centered ($\E(X_{i})=0$ for all $i$),
real-valued bounded random variables, which are $\rho$-mixing in the sense that
\[
\rho(k)=\sup_{j \geq 1} \sup_{v >u \geq j+k}\rho\big (\sigma(X_{i},1\leq
i\leq j),\sigma(X_{u}, X_{v} )\big )\rightarrow0\,,\,\text{as
$k\rightarrow\infty$}\,,
\]
where $\sigma(X_{t},t\in A)$ is the $\sigma$-field generated by the r.v.'s
$X_{t}$ with indices in $A$ and we recall that the maximal correlation
coefficient $\rho({\mathcal{U}},{\mathcal{V}})$ between two $\sigma$-algebras
is defined by
\[
\rho(\mathcal{U},\mathcal{V})=\sup\{|\mathrm{corr}(X,Y)|:X\in{\mathbb{L}}%
^{2}(\mathcal{U}),Y\in{\mathbb{L}}^{2}(\mathcal{V})\}\,.
\]
In this section we shall also assume that the r.v.'s $(X_i)_{i  \geq 1}$ satisfies the following set of assumptions 
\[
(H):=%
\begin{cases}
1) $ $   \Theta = \sum_{k \geq 1} k \rho(k) < \infty \, .\\
2) $  $  \mbox{For any $n \geq 1$, }   \displaystyle  C_n := \max_{1 \leq \ell \leq n} \frac{ \sum_{i=\ell}^n \E (X_i^2)}{\E (S_n -S_{\ell-1})^2}  < \infty\, . 
\end{cases}
\]

\begin{Remark} Note that in $(H_2)$ necessarily $C_n \geq 1$.  In many cases of interest the sequence $(C_n)_n$ is bounded: for example, when $X_i=f_i(Y_i)$ where $Y_i$ is a Markov chain satisfying $\rho_Y(1) < 1$, then according to \cite[Proposition 13]{P12},   $  C_n \leq (1+ \rho_Y(1) ) ( 1 - \rho_Y(1) )^{-1}$. Here  $(\rho_Y(k))_{k \geq 0}$ is the sequence of $\rho$-mixing coefficients of the Markov chain $(Y_i)_i$. 
\end{Remark}

\begin{Corollary} \label{corappli1}   Let $(X_i)_{i  \geq 1}$ be a sequence of centered bounded real-valued random variables such that $(H)$ is satisfied. Let $V_n = {\rm Var} (S_n)$ and $K_n = \max_{1 \leq i \leq n} \Vert X_i \Vert_{\infty}$. Then for any positive integer $n$, 
\[
W_1(P_{S_{n}}, G_{V_n})  \ll K_n ( 1+  C_n \log (1+ C_n V_n) )  \, .
\]
\end{Corollary}
\begin{Remark} If the sequences $(C_n)_n$ and $(K_n)_n$ are bounded and $V_n \rightarrow \infty$, then Corollary \ref{corappli1} provides a rate in the central limit theorem for $S_n / \sqrt{V_n}$. More precisely, 
\[
W_1(P_{S_{n}/\sqrt{V_n}}, G_{1})  = O ( V_n^{-1/2}   \log ( V_n)  ) \  \text{ and } \   \Vert F_n  - \Phi   \Vert_{\infty} = O (  V_n^{-1/4}  \sqrt{\log (V_n) } ) \, .
\]
where $F_n$ is the c.d.f. of $S_n /{\sqrt V_n}$ (the second inequality follows from \eqref{ineBEW1}).  Note that   the above upper bounds hold even if we do not require a linear growth of the variance $V_n $ as it is imposed for instance in 
\cite[Theorem 3.1]{Wang-Hu} and of course, in the stationary case, in \cite{Zu,Rio95,Ti80}.
\end{Remark}

\subsection{Sequential dynamical systems} \label{subsectionseqdyn}
The term sequential dynamical system, introduced by Berend and Bergelson \cite{BB}, refers to a non-stationary system defined by the composition of deterministic maps $T_k \circ T_{k-1} \circ \cdots \circ T_1$ acting on a space $X$. 

More precisely, we consider here the setting described by Conze and Raugi \cite{CR} and Haydn et al. \cite{HNTV}. 
Let $(T_k)_{k \geq 1}$ be a sequence of maps from $X$ to $X$, where $X$ is either a compact subset of ${\mathbb R}^d$ or the $d$-dimensional torus ${\mathbb T}^d$. Let also $m$ be the Lebesgue measure defined on the Borel $\sigma$-algebra ${\mathcal B}$ of $X$, normalized in such a way that $m(X)=1$. We assume that each $T_k$ is non singular with respect to $m$ i.e. $m(A) >0 \Longrightarrow m(T(A))>0$. 

Let $P_k$ be the Perron-Frobenius operator, that is the adjoint of the composition by $T_k$: for any $f \in {\mathbb L}_1(m), g \in {\mathbb L}_\infty(m)$, 
$$
  \int_X f(x) \, g \circ T_k(x) \, m(dx)= \int_X (P_k f)(x)  \,  g(x) \, m(dx) \, .
$$ 
Let also $\tau_k= T_k \circ T_{k-1} \circ \ldots \circ T_1$ and $ \pi_k = P_k \circ P_{k-1}  \circ \ldots \circ P_1$, and note that $\pi_k$ is the Perron-Frobenius operator of $\tau_k$. 

Let ${\mathcal V} \subset {\mathbb L}_\infty(m)$, $(1 \in {\mathcal V})$,  be a Banach space of functions from $X$ to ${\mathbb R}$  with norm $\| \cdot \|_v$, such that $\|\phi \|_\infty \leq \kappa_1 \|\phi \|_v$ for some $\kappa_1 >0$. We assume moreover that if $\phi_1, \phi_2$ are two functions in ${\mathcal V}$, then the usual product $\phi_1 \phi_2$ belongs to ${\mathcal V}$ and satisfies $\|\phi_1 \phi_2 \|_v \leq \kappa_2 \|\phi_1\|_v \|\phi_2\|_v $ for some $\kappa_2 >0$. In what follows, we set $\kappa = \max(\kappa_1, \kappa_2)$. Typical examples of Banach spaces ${\mathcal V}$ are the space $BV$ of functions with bounded variation on a compact interval of ${\mathbb R}$, or the space  ${\mathcal H}_\alpha$ of $\alpha$-H\"older function on a compact set of ${\mathbb R}^d$, equipped with their usual norms. 

We now recall the properties (DEC) and (MIN) introduced in \cite{CR} (we use the formulation of 
\cite{HNTV}): 

\medskip

\noindent {\bf Property (DEC):}  There exist two constants $C>0$ and $\gamma \in (0,1)$ such that: for any positive integer $n$, any $n$-tuple $(j_1, \ldots , j_n)$ of positive integers, and any $f \in {\mathcal V}$,
\[
   \left \|  P_{j_n} \circ \cdots \circ  P_{j_1} (f -m(f))   \right \|_v  \leq C \gamma^n \|f- m(f)\|_v \, .
\]

\medskip

\noindent {\bf Property (MIN):} There exist $\delta >0$ and $\gamma \in (0,1)$ such that: for any positive integer $n$,  and any $n$-tuple $(j_1, \ldots , j_n)$ of positive integers, we have the uniform lower bound
\[
    \inf_{x \in X} P_{j_n} \circ \cdots \circ P_{j_1} 1 (x) \geq \delta \, .
\]

\medskip

The main result of this subsection is the following corollary.

\begin{Corollary} \label{corsequential}
Let $(\phi_n)_{n \geq 1}$ be a sequence of functions in ${\mathcal V}$ such that 
$\sup_{n\geq 1} \|\phi_n \|_v < \infty$. Let 
\[
S_n= \sum_{k=1}^n \left (\phi_k (\tau_k) - m(\phi_k (\tau_k)) \right) \, , \quad and \quad 
V_n= \int_{X} S_n^2 (x) \, m(dx) \, .
\]
Assume that the properties  (DEC) and (MIN)  are satisfied. Then, on the probability  space 
$(X, {\mathcal B}, m)$,
\[
W_1(P_{S_n}, G_{V_n}) \ll \log (n+1) \log(2+ V_n)  \, . 
\]
\end{Corollary}
\begin{Remark} Under the assumptions of Corollary \ref{corsequential}, we derive that 
\[
W_1(P_{S_n/ {\sqrt{V_n}}}, G_{1}) \ll  V_n^{-1/2}\log (n+1) \log(2+ V_n) 
\]
and 
\[
\Vert F_n - \Phi \Vert_{\infty} \ll  \Big ( V_n^{-1/2}\log (n+1) \log(2+ V_n)    \Big )^{1/2} \, , 
\]
where $F_n$ is the cdf of $S_n/ {\sqrt V_n}$ (the second inequality follows from \eqref{ineBEW1}). In particular, 
Corollary \ref{corsequential} provides a rate in the central limit theorem for $S_n / \sqrt V_n$ as soon as $(\log n \log \log n) / \sqrt V_n  \rightarrow 0$ as $n \rightarrow \infty$.
\end{Remark}

\section{Proofs} \label{sectionproofs}
\setcounter{equation}{0}

\subsection{Proof of Theorem \ref{Th1first}}

The proof is based on the following proposition:

\begin{Proposition} \label{Mainprop}
Let  $\delta$ be a positive real and  denote by 
 $t_{\ell,n} =\big (  V_n - V_{\ell } + \delta^2 \big )^{1/2}$.  Let $p \in ]2,3]$ and $r \in (0, p]$.  Then, 
 there exist  positive constants $C_{r,p}$ depending on $(r,p)$ and $\kappa_r$ depending on $r$ such that for every positive integer $n$,
\begin{multline} \label{boundzetaprop}
\zeta_r(P_{M_{n}}, G_{V_n})  \leq   4\sqrt{2}  \delta^r + C_{r,p}   \Bigl  \{   \sum_{k=1}^n    \Big (    \frac{1}{ t_{k,n}^{3-r}}  \E \big (  \xi_k^2  \min ( \kappa_r  t_{k,n}  ,  |\xi_k|  ) \big ) + \frac{\sigma_k^4}{t_{k,n}^{4-r}}   \Big )   +   \sum_{ \ell =2}^n \frac{U_{\ell,n} (p) }{ (t_{\ell-1,n} )^{p-r}}  \Bigr \}  
 \, ,
\end{multline}
where, for $\ell \geq 2$,  $
 U_{\ell,n} (p) $ is defined in \eqref{defU}. 
\end{Proposition}
\begin{Remark} \label{Rollin} When $r=1$, $p=3$ and  $U_{\ell,n} (p)=0$ for any $\ell$, our bound is  similar to the one stated in \cite[Theorem 2.1]{Rollin18}. However our quantity $ \sum_{ \ell =2}^n (t_{\ell-1,n} )^{r-p} U_{\ell,n} (p)  $ can be handled in many cases (see Section \ref{sectionapplications}) while his condition $V_n^{-1} \langle M \rangle_n =1$ a.s. is very restrictive. 
\end{Remark}
We end the proof of the theorem with the help of this proposition taking  $\delta= a \delta_n$. Hence we shall give an upper bound for 
\[
\sum_{k=1}^n    \Big (    \frac{1}{ t_{k,n}^{3-r}}  \E \big (  \xi_k^2  \min ( \kappa_r  t_{k,n}  ,  |\xi_k|  ) \big ) + \frac{\sigma_k^4}{t_{k,n}^{4-r}}   \Big ) \, , 
\]
where  $t_{k,n} = (  a^2\delta_n^2 + \sigma^2_{k+1} + \cdots + \sigma^2_{n} )^{1/2} $. 
With this aim note first that 
\[
 \frac{1}{ t_{k,n}^{3-r}}  \E \big (  \xi_k^2  \min ( \kappa_r  t_{k,n}  ,  |\xi_k|  ) \big ) \leq  \frac{\sigma_k^2}{ t_{k,n}^{3-r}}  \psi_n ( \kappa_r \delta_n^{-1}t_{k,n} ) \, ,
\]
where $\psi_n(t)$ is defined in \eqref{defpsi}. 
Let ${\tilde \sigma}_k = { \sigma}_k/ \delta_n$. Note that since  ${\tilde \sigma}_k \leq 1$,
\[
\frac{\sigma_k^2}{ t_{k,n}^{2}}  = \frac{{\tilde \sigma}_k^2}{ a^2 +{\tilde \sigma}^2_{k+1} + \cdots + {\tilde \sigma}^2_{n} }  \leq \frac{ \alpha {\tilde \sigma}_k^2}{ a^2 +{\tilde \sigma}_k^2 + \alpha \sum_{\ell = k+1}^n {\tilde \sigma}^2_{\ell} }  \, ,
\]
where $\alpha = (a^2+1)/a^2$.  Let $u_k = a^2 + \alpha \sum_{\ell = k+1}^n {\tilde \sigma}^2_{\ell} $. It follows that 
\[
\frac{\sigma_k^2}{ t_{k,n}^{2}}  \leq \frac{u_{k-1} - u_k}{ (u_{k-1} - u_k)/\alpha +u_k } =  \frac{\alpha ( u_{k-1} - u_k) }{ (u_{k-1} - u_k ) + \alpha u_k } =   \frac{\alpha a_k }{ a_k + \alpha } 
\]
where 
\[
a_k = ( u_{k-1} - u_k) /u_k \, .
\]
But since $a^2 \geq 1$ we have $\alpha \leq 2$. Hence, for any $x \geq 0$, 
\[
\frac{\alpha x}{x+\alpha} \leq \log (1+x) \, ,
\]
implying that
\beq \label{trivial1}
\frac{\sigma_k^2}{ t_{k,n}^{2}} \leq  \log ( 1+a_k) =  \log (u_{k-1} /u_k) \, .
\eeq
It follows that, if $r \geq 1$, since $t \mapsto \psi_n(t)$ is non decreasing and $t_{k,n}^2 \leq \delta_n^2 u_k$ (since $\alpha \geq 1$), \begin{multline*}
 \frac{\sigma_k^2}{ t_{k,n}^{3-r}} \psi_n ( \kappa_r \delta_n^{-1}t_{k,n} )  =  \frac{\sigma_k^2}{ t_{k,n}^{2}} \psi_n ( \kappa_r \delta_n^{-1}t_{k,n} )   t_{k,n}^{r-1} \leq    2  \log (\sqrt{u_{k-1}} / \sqrt{u_k}) \psi_n ( \kappa_r  \sqrt{u_{k}} )  \delta_n^{r-1} u_{k}^{(r-1)/2} \\
 \leq 2  \psi_n ( \kappa_r \sqrt{u_{k}} )  \delta_n^{r-1} u_{k}^{(r-1)/2} \int_{\sqrt{u_k}}^{\sqrt{u_{k-1}}} \frac{1}{x} dx  \leq 2 \delta_n^{r-1} \int_{\sqrt{u_k}}^{\sqrt{u_{k-1}}} \frac{\psi_n ( \kappa_r  x ) }{x^{2-r}} dx \, .
\end{multline*}
Hence, if $r \geq 1$,
\begin{align} \label{trivial2}
\sum_{k=1}^n  \frac{\sigma_k^2}{ t_{k,n}^{3-r}} \psi_n ( \kappa_r \delta_n^{-1}t_{k,n} )   &  \leq 2  \delta_n^{r-1} \int_{a}^{ \sqrt{a^2 + \alpha \sum_{\ell =1}^n {\tilde \sigma}^2_{\ell}} }\frac{\psi_n ( \kappa_r  x ) }{x^{2-r}}dx \nonumber  \\
& \leq  2 \delta_n^{r-1}   \int_{a}^{\sqrt{ v_n (a)  / \delta_n^2}}  \frac{\psi_n ( \kappa_r  x ) }{x^{2-r}}dx \, .
\end{align}
We study now the case $r <1$. With this aim, note first that taking into account that  ${\tilde \sigma}^2_{k}  \leq 1$, $\alpha \leq 2$ and that $a \geq 1$, we have
\begin{equation} \label{secondboundtk}
t_{k,n}^2 = \delta_n^2  \Big ( a^2 + \sum_{\ell = k+1}^n {\tilde \sigma}^2_{\ell} \Big )  \geq    a^2 ( a^2 + \alpha)^{-1}  \delta_n^2   u_{k-1} \geq     \delta_n^2   u_{k-1}   /3 \, ,
\end{equation}
(for the first inequality, use the fact that $a^2 ( a^2 + \alpha)^{-1}  \leq \alpha^{-1}$). 
When $r <1$,  taking into account the upper bound \eqref{secondboundtk}, we then derive 
\[
 \frac{\sigma_k^2}{ t_{k,n}^{3-r}} \psi_n ( \kappa_r \delta_n^{-1}t_{k,n} )   \leq    2 \times 3 ^{(1-r)/2}   \delta_n^{r-1} u_{k-1}^{(r-1)/2} \psi_n ( \kappa_r \sqrt{u_{k}} )   \log (\sqrt{u_{k-1}} / \sqrt{u_k})  \, .
 \]
Hence, when $r <1$, 
 \[
\sum_{k=1}^n  \frac{\sigma_k^2}{ t_{k,n}^{3-r}} \psi_n ( \kappa_r \delta_n^{-1}t_{k,n} )  \leq   2 \times 3 ^{(1-r)/2}  \delta_n^{r-1}   \int_{a}^{\sqrt{ v_n (a) / \delta_n^2}}  \frac{\psi_n ( \kappa_r  x ) }{x^{2-r}}dx \, .
\]

The bound \eqref{secondboundtk} and \eqref{trivial1} also implies that, for any $r \leq 2$, 
 \begin{multline*}
\sum_{k=1}^n  \frac{\sigma_k^4}{ t_{k,n}^{4-r}}  \leq \delta_n^2  \sum_{k=1}^n  \frac{\sigma_k^2}{t_{k,n}^2}  \times  \frac{1}{ t_{k,n}^{2-r}}  \leq 3^{(2-r)/2} \delta_n^r \sum_{k=1}^n  \frac{\sigma_k^2}{t_{k,n}^2}  \times  \frac{1}{u_{k-1}^{(2-r)/2}}  
\\  \leq 2 \times 3 ^{(2-r)/2}\delta_n^r \sum_{k=1}^n  \log (\sqrt{u_{k-1}} / \sqrt{u_k})  \times  \frac{1}{u_{k-1}^{(2-r)/2}}   =  2 \times 3 ^{(2-r)/2}\delta_n^r \sum_{k=1}^n  \frac{1}{u_{k-1}^{(2-r)/2}}   \int_{\sqrt{u_k}}^{\sqrt{u_{k-1}}} \frac{1}{x} dx \\
\leq 2 \times 3 ^{(2-r)/2} \delta_n^r \sum_{k=1}^n   \int_{\sqrt{u_k}}^{\sqrt{u_{k-1}}} \frac{1}{x^{3-r}} dx \leq  2 \times 3 ^{(2-r)/2} \delta_n^r    \int_{a}^{\sqrt{u_{0}}} \frac{1}{x^{3-r}} dx  \, .
\end{multline*}
When $r >2$, we use the fact that $t_{k,n}^2 \leq \delta_n^2 u_k$ to derive that 
\[
\sum_{k=1}^n  \frac{\sigma_k^4}{ t_{k,n}^{4-r}}  \leq  2   \delta_n^r    \int_{a}^{\sqrt{u_{0}}} \frac{1}{x^{3-r}} dx \, .
\]
All these considerations end the proof of Theorem \ref{Th1first}. It remains to prove Proposition \ref{Mainprop}.

\medskip

\noindent {\it Proof of Proposition \ref{Mainprop}.}  Let $\Y$ be a sequence of ${\mathcal N}(0, \sigma_i^2)$-distributed independent random
variables, independent of the sequence $(\xi_i)_{i \in {\mathbb N}}$. For $n>0$, let
$T_{n} = \sum_{j=1}^n Y_j$. Let also $Z$ be a ${\mathcal N}(0, \delta^2)$-distributed random variable
independent of $(\xi_i)_{i \in {\mathbb N}}$ and $(Y_i)_{i \in
{\mathbb N}}$. Using Lemma 5.1 in \cite{DMR} together with the fact that, for any real $c$, 
$\zeta_r(P_{cX}, P_{cY})= |c|^r \zeta_r(P_{X}, P_{Y})$, we derive
that for any $r$ in $]0, p]$, \begin{equation} \label{conv} \zeta_r
( P_{M_{n}} , P_{T_{n}} ) \leq 2 \zeta_r ( P_{M_{n}}
* P_{Z} , P_{T_{n}}
* P_{Z} ) + 4\sqrt{2}   \delta^r \, .\end{equation}
Consequently it remains to bound up
$$
\zeta_r ( P_{M_{n}}
* P_{Z} , P_{T_{n}}
* P_{Z} ) = \sup_{f \in {\Lambda_r}} {\mathbb E} ( f (M_{n}
+ Z) - f (T_{n} + Z)) \, .
$$Recall that $V_n = \sum_{i=1}^n \sigma_i^2$ and, for any $k\leq n$, set
$$
f_{V_n-V_k} (x) =  {\mathbb E} ( f (x+ T_{n}-T_{k} + Z)  )  .
$$
Then, from the independence of the above sequences, \[ \bkE ( f ( M_{n} + Z ) - f ( T_{n} + Z ) )
= \sum_{ k=1}^{n} D_{k} \, , \]where \[ D_{k} = \bkE  \big(
f_{V_n-V_k}(M_{k-1}   +  \xi_k ) -f_{V_n-V_k}(M_{k-1}  +
Y_k ) \big) \, .
\]
By the Taylor formula, we get 
\begin{multline*}
f_{V_n-V_k}(M_{k-1}   +  \xi_k ) -f_{V_n-V_k}(M_{k-1}  +
Y_k )  \\= f'_{V_n-V_k} (M_{k-1}) ( \xi_k- Y_k ) + \frac{1}{2}  f''_{V_n-V_k}(M_{k-1})  ( \xi^2_k- Y^2_k ) - \frac{1}{6}  f^{(3)}_{V_n-V_k} (M_{k-1}) (Y^3_k ) + R_k \, , 
\end{multline*}
where
\[
R_k \leq  \xi_k^2 \Big ( \Vert  f''_{V_n-V_k} \Vert_{\infty} \wedge \frac{1}{6} \Vert  f^{(3)}_{V_n-V_k} \Vert_{\infty}  |\xi_k| \Big ) + \frac{1}{24} \Vert  f^{(4)}_{V_n-V_k} \Vert_{\infty}  Y_k^4  \, .
\]
Using the fact that $( \xi_k)_{k \in {\mathbb N}}$ is a sequence of martingale differences independent of the sequence of iid Gaussian random variables $(Y_k)_{k \in {\mathbb N}}$, we then get
\beq
\label{telessum} \bkE ( f ( M_{n} + Z ) - f ( T_{n} + Z ) )
= \frac{1}{2}   \sum_{ k=1}^{n} \E \big (  f''_{V_n-V_k}(M_{k-1})  ( \xi^2_k- Y^2_k )  \big ) + \sum_{ k=1}^{n}   \E (R_k)  \, . \eeq
Note first that 
\[
\E (R_k) \leq    \E \Big ( \xi_k^2 \Big ( \Vert  f''_{V_n-V_k} \Vert_{\infty} \wedge \frac{1}{6} \Vert  f^{(3)}_{V_n-V_k} \Vert_{\infty}  |\xi_k| \Big )  \Big ) + \frac{\sigma_k^4}{8} \Vert  f^{(4)}_{V_n-V_k} \Vert_{\infty}   \, .
\]
Recall the notation $t_{k,n} = (  \delta^2 + \sigma^2_{k+1} + \cdots + \sigma^2_{n} )^{1/2} $. By Lemma 6.1 in \cite{DMR}, we have that for any integer $i \geq 1$, 
\beq \label{lma1DMR}
\Vert  f^{(i)}_{V_n-V_k}   \Vert_{\infty} \leq c_{r,i} t_{k,n}^{r-i} \, .
\eeq
Hence, setting $\kappa_r =  6 c_{r,2} / c_{r,3} $, we get 
\beq \label{b1Rkprop1}
\E (R_k) \leq  \frac{c_{r,3}}{6} \times \frac{1}{ t_{k,n}^{3-r}}  \E \Big (  \xi_k^2  \min ( \kappa_r  t_{k,n}  ,  |\xi_k|  ) \Big ) + \frac{c_{r,4}}{8} \frac{\sigma_k^4}{t_{k,n}^{4-r}}  \, . 
\eeq
For $r=1$, we can take  $\kappa_r =6$, $c_{r,3}=1$ and $c_{r,4} =8/5$. 

\medskip

We study now the quantity 
$ \sum_{ k=1}^{n}  \E \big ( f''_{V_n-V_k}(M_{k-1})  ( \xi^2_k- Y^2_k ) \big )$. 
With this aim let us consider a sequence $(Y_k')$ of real-valued random variables  independent of $(Y_k) $ and $(\xi_k)$ and such that ${ \mathcal L} (Y_k') ={ \mathcal L} (Y_k) $. Note first that   
\[
 \E   \big (    ( f''_{V_n-V_k}(M_{k-1} + Y_k' ) -  f''_{V_n-V_k}(M_{k-1} )  )  ( \xi^2_k- Y^2_k )  \big )  =  \E   \big (     f^{(3)}_{V_n-V_k}(M_{k-1}  )  Y_k'    ( \xi^2_k- Y^2_k )  \big )   + \E (R_k')   \, ,
\]
where, by taking into account    \eqref{lma1DMR} and the independence between $(Y_k')_k$ and $(\xi_k, Y_k)_k$, 
\[
 \E (|R_k'|) ) \leq       \Vert  f^{(4)}_{V_n-V_k} \Vert_{\infty}  \E \big | (Y_k')^2   ( \xi^2_k- Y^2_k )  \big |    
\leq  2 c_{r,4}  \frac{   \sigma_k^4 }{t_{k,n}^{4-r}}    \, . 
\]
Since $\E (Y_k') =0$ and  $(Y_k')_k$ is independent of  $(\xi_k, Y_k)_k$, we get 
\beq \label{b2prop1}
\sum_{ k=1}^{n}  \Big | \E   \big (    ( f''_{V_n-V_k}(M_{k-1} + Y_k' ) -  f''_{V_n-V_k}(M_{k-1} )  )  ( \xi^2_k- Y^2_k )  \big )  \Big |  \leq  2 c_{r,4} \sum_{ k=1}^{n}   \frac{   \sigma_k^4 }{t_{k,n}^{4-r}} \, .
\eeq
Now
 \begin{multline*}
 \E   \big (     f''_{V_n-V_k}(M_{k-1} + Y_k' )  ( \xi^2_k- Y^2_k )  \big )  =   \E   \big (     f''_{V_n-V_{k-1}}(M_{k-1}  )  ( \xi^2_k- Y^2_k )  \big )   \\
  = \sum_{\ell =2}^k   \E   \Big (   \big (     f''_{V_n-V_{k-1}}(M_{ \ell-1}  + T_{k-1} - T_{\ell -1} ) -   f''_{V_n-V_{k-1}}(M_{ \ell-2}  + T_{k-1} - T_{\ell -2} )     \big ) ( \xi^2_k- Y^2_k )  \Big )  \\
    = \sum_{\ell =2}^k   \E   \Big (   \big (     f''_{V_n-V_{\ell-1}}(M_{ \ell-1}   ) -   f''_{V_n-V_{\ell-1}}(M_{ \ell-2}  + T_{\ell-1} - T_{\ell -2} )     \big ) ( \xi^2_k- Y^2_k )  \Big ) \, . 
\end{multline*}
Hence, by using  Lemma 6.1 in \cite{DMR}, there exists a positive constant $c_{r,p}$ depending on $(r,p)$ such that for any $n \geq 1$, 
\begin{multline}  \label{b3prop1}
\sum_{k=1}^n  \E   \big (     f''_{V_n-V_k}(M_{k-1} + Y_k' )  ( \xi^2_k- Y^2_k )  \big )  \\
=   \sum_{\ell=2}^n    \E   \Big (   \big (     f''_{V_n-V_{\ell-1}}(M_{ \ell-1}   ) -   f''_{V_n-V_{\ell-1}}(M_{ \ell-2}  + T_{\ell-1} - T_{\ell -2} )     \big ) \sum_{k=\ell}^n  ( \E_{\ell -1} (\xi^2_k)- \sigma^2_k )  \Big )  \\
\leq c_{r,p}  \sum_{\ell=2}^n    \frac{1}{(V_n -V_{\ell -1} +  \delta^2 )^{(p-r)/2} } \Big  \Vert  | \xi_{\ell -1} - Y_{\ell-1}|^{p-2} \Big |  \sum_{k=\ell}^n  ( \E_{\ell -1} (\xi^2_k)- \sigma^2_k )   \Big | \Big  \Vert_1 \, .
\end{multline}
Starting from \eqref{telessum} and taking into account the upper bounds \eqref{b1Rkprop1},  \eqref{b2prop1} and \eqref{b3prop1}, the desired inequality follows since for any integer $\ell \in [2,n]$ and any $p\in [2,3]$,  we have 
$\E ( | Y_{\ell-1}| ^{p-2} ) \leq ( \E | Y_{\ell-1}|) ^{p-2}  \leq \sigma_{\ell-1}^{p-2} $. \qed 

\subsection{Proof of Corollary \ref{corapplilinearcondonspectral}} For any $k \geq 1$, let ${\mathcal F}_k = \sigma (X_1, \ldots, X_k)$ and ${\mathcal F}_0= \{  \emptyset, \Omega\}$. Write first
\[
S_n = \sum_{k=1}^n  ( \E_k (S_n) - \E_{k-1} (S_n) ) =: \sum_{k=1}^n d_{k,n} \, .
\]
Note that $(d_{k,n})_{1 \leq k \leq n}$ is a triangular array of martingale differences with respect to $({\mathcal F}_k )_{k \geq 1}$ and that $V_n= \sum_{k=1}^n \E (d_{k,n}^2)  = \E (S_n^2)$. Hence,  setting $\delta_n = \max_{1 \leq k \leq n}
\Vert  d_{k,n}\Vert_2$ and applying Proposition \ref{Mainprop} we get that, for any $a \geq 1$, 
\beq \label{startlin}
W_1(P_{S_{n}}, G_{V_n})  \ll  a \delta_n +      \sum_{k=1 }^n  \Big (  \frac{ \E( |d_{k,n}|^p)}{ B^{(p-1)/2}_{ k+1 , n} (a)}  +  \frac{\sigma_{k,n}^4}{B^{3/2}_{ k+1 , n} (a)}   \Big ) + \sum_{ \ell =2}^n \frac{1}{ B^{(p-1)/2}_{ \ell , n} (a)  } U_{\ell,n} (p) 
 \, ,
\eeq
where $\sigma_{k,n} = \Vert d_{k,n} \Vert_2$ and 
\[
B_{ \ell , n} (a) =\sum_{k=\ell }^n \E (d_{k,n}^2) +  a^2 \delta_n^2 \,  \mbox{ and } \,  U_{\ell,n} (p)  = \Big \Vert( | d_{\ell -1,n} | \vee   \sigma_{\ell -1,n} ) \Big |^{p-2} \sum_{k = \ell}^n ( \E_{\ell-1} ( d_{k,n}^2) - \sigma_{k,n}^2) \Big |   \Big \Vert_1
\, .
\]
Proceeding as in the proof of Theorem \ref{Th1first}, we get that 
\beq \label{boundterme4}
\sum_{k=1 }^n \frac{\sigma_{k,n}^4}{B^{3/2}_{ k+1 , n} (a)} \ll \delta_n \, .
 \eeq 
  Next, 
setting $\alpha = (a^2+1)/a^2$,  note that 
\[
B_{ k+1 , n} (a) \geq \alpha^{-1} \Big ( a^2 \delta^2_n  + \sigma_{k,n}^2 + \alpha  \sum_{\ell =k+1}^n  \sigma_{\ell,n}^2 \Big ) 
\geq 2^{-1} B_{ k , n} (a)  \, .
\]
Note also that 
\[
U_{\ell,n} (p)  \leq 2   \Vert d_{\ell -1,n} \Vert_p^{p-2} \big  \Vert \sum_{k = \ell}^n ( \E_{\ell-1} ( d_{k,n}^2) - \sigma_{k,n}^2) \big  \Vert_{p/2} \, .
\]
But, setting $A_{k,n} = \E_{k} (S_n- S_k) $, note that the following decomposition is valid:
\begin{equation} \label{decmart1}
d_{k,n } = X_{k,n} + A_{k,n} - A_{k-1,n} \, .
\end{equation}
Hence
\[
3^{1-p} \Vert d_{\ell ,n} \Vert^p_p \leq     \Vert X_{\ell ,n} \Vert^p_p  + \Vert A_{\ell,n} \Vert^p_p  + \Vert A_{\ell-1,n} \Vert^p_p \ \leq   | \alpha_{n , \ell } |^p \Vert Y_{0} \Vert^p_p  + 2 \Big (  \sum_{i = \ell}^n  | \alpha_{i,n}  |   \Vert \E (Y_i | {\mathcal G}_{\ell-1} ) \Vert_p  \Big )^p   \, .
\]
But, by convexity, setting $\beta_i = \Vert \E ( Y_i | {\mathcal G}_{\ell-1}  )\Vert_p \big (   \sum_{u=\ell}^n \Vert \E ( Y_u | {\mathcal G}_{\ell-1}  )\Vert_p \big )^{-1}$, we get 
\begin{multline*}
 \Big (  \sum_{i=\ell}^n | \alpha_{i,n}  | \Vert \E ( Y_i | {\mathcal G}_{\ell-1}  )\Vert_p\Big )^p  \leq  \sum_{i=\ell}^n |\alpha_{i,n} |^p  \beta_i^{1-p} \Vert \E ( Y_i | {\mathcal G}_{\ell-1}  )\Vert^p_p \\
\leq \Big (   \sum_{u=1}^n \Vert \E ( Y_u | {\mathcal G}_{0}  )\Vert_p  \Big )^{p-1}  \sum_{i=\ell}^n |\alpha_{i,n} |^p \Vert \E ( Y_i | {\mathcal G}_{\ell-1}  )\Vert_p   \, ,
\end{multline*}
implying that 
\beq  \label{b1dln}
\Vert d_{\ell ,n} \Vert^p_p \ll \Big (   \sum_{u=0}^n \Vert \E ( Y_u | {\mathcal G}_{0}  )\Vert_p  \Big )^{p-1}    \sum_{i=\ell}^n |\alpha_{i,n} |^p \Vert \E ( Y_i | {\mathcal G}_{\ell}  )\Vert_p \, . 
\eeq
It follows  that 
\beq  \label{b2dln}
\max_{1 \leq \ell \leq n}\Vert d_{\ell ,n} \Vert^{p-2}_p \ll \max_{1 \leq i \leq n} |\alpha_{i,n} |^{p-2}\Big (   \sum_{u=0}^n \Vert \E ( Y_u | {\mathcal G}_{0}  )\Vert_p  \Big )^{p-2} :=  \max_{1 \leq i \leq n} |\alpha_{i,n} |^{p-2} \eta_n^{p-2} \, . 
\eeq
On another hand
\begin{multline} \label{p1martpart}
\Big \Vert \sum_{k = \ell}^n ( \E_{\ell-1} ( d_{k,n}^2) - \sigma_{k,n}^2)   \Big \Vert_{p/2} = \Big  \Vert  \E_{\ell-1} \Big (   \sum_{k = \ell}^n d_{k,n} \Big )^2 -  \E \Big (   \sum_{k = \ell}^n d_{k,n} \Big )^2   \Big \Vert_{p/2} \\
=  \Vert  \E_{\ell-1} ( S_n - \E_{\ell -1} (S_n) )^2 -   \E ( S_n - \E_{\ell -1} (S_n) )^2 \Vert_{p/2}  \\
  \leq   \Vert  \E_{\ell-1} ( S_n  -S_{\ell -1} )^2 -\E ( S_n  -S_{\ell -1} )^2 \Vert_{p/2}  + 2  \Vert  \E_{\ell-1} ( S_n  -S_{\ell -1} )  \Vert^2_{p}  \, .
\end{multline}
Note that 
\begin{multline}  \label{p2martpart}
  \Vert  \E_{\ell-1} ( S_n  -S_{\ell -1} )^2 -\E ( S_n  -S_{\ell -1} )^2 \Vert_{p/2}   \leq 2 \sum_{i= \ell}^n \sum_{j= i}^n   \Vert  \E_{\ell-1} ( X_{i,n} X_{j,n} ) - \E ( X_{i,n} X_{j,n} ) \Vert_{p/2} \\
\leq  2 \sum_{i= \ell}^n \sum_{j= i}^n  \alpha_{i,n} \alpha_{j,n}   \Vert  \E ( Y_{i} Y_{j} | {\mathcal G}_{\ell -1} ) - \E ( Y_{i} Y_{j} ) \Vert_{p/2} \\
\leq   2 \sum_{i= \ell}^n \sum_{j= i}^{2i-\ell}  \alpha_{i,n} \alpha_{j,n}   \Vert  \E ( Y_{i} Y_{j} | {\mathcal G}_{\ell -1} ) - \E ( Y_{i} Y_{j} ) \Vert_{p/2}  +   4 \sum_{i= \ell}^n \sum_{j= 2i-\ell +1}^n   \alpha_{i,n} \alpha_{j,n}   \Vert Y_{i}  \E ( Y_{j} | {\mathcal G}_{i} )  \Vert_{p/2}  \, .
\end{multline}
Hence by stationarity,
\begin{multline*}
 \Vert  \E_{\ell-1} ( S_n  -S_{\ell -1} )^2 -\E ( S_n  -S_{\ell -1} )^2 \Vert_{p/2}  \\
\leq   4 \Big (  \sum_{i= \ell}^n \sum_{j= i}^{n \wedge (2i-\ell)}  \alpha_{i,n} \alpha_{j,n}   \lambda_{i-\ell +1} +    \sum_{i= \ell}^n \sum_{j= 2i-\ell +1}^n   \alpha_{i,n} \alpha_{j,n}    \lambda_{j-i}  \Big )  \, .
\end{multline*}
It follows that 
\begin{multline*}
  \Vert  \E_{\ell-1} ( S_n  -S_{\ell -1} )^2 -\E ( S_n  -S_{\ell -1} )^2 \Vert_{p/2}  \\
\leq   2  \Big ( \sum_{i= \ell}^n  \alpha^2_{i,n}  (i-\ell +1)   \lambda_{i-\ell +1} +  2  \sum_{j= \ell}^{n }  \alpha^2_{n,j}   \sum_{u= [(j-\ell)/2]}^{j - \ell} \lambda_{u}   +  \sum_{i= \ell}^n   \alpha^2_{i,n} \sum_{u= i-\ell +1}^{n-i}    \lambda_{u} \Big )  \, .
\end{multline*}
In addition, setting $\beta_i = \Vert \E ( Y_i | {\mathcal G}_{\ell-1}  )\Vert_p \big (   \sum_{u=\ell}^n \Vert \E ( Y_u | {\mathcal G}_{\ell-1}  )\Vert_p \big )^{-1}$, we get by convexity,  
\begin{multline}  \label{p3martpart}
  \Vert  \E_{\ell-1} ( S_n  -S_{\ell -1} )  \Vert^2_{p}  = \Big (  \sum_{i=\ell}^n \alpha_{i,n} \Vert \E ( Y_i | {\mathcal G}_{\ell-1}  )\Vert_p\Big )^2  \leq  \sum_{i=\ell}^n \alpha^2_{i,n}  \beta_i^{-1} \Vert \E ( Y_i | {\mathcal G}_{\ell-1}  )\Vert^2_p \\
\leq  \sum_{u=1}^n \Vert \E ( Y_u | {\mathcal G}_{0}  )\Vert_p   \sum_{i=\ell}^n \alpha^2_{i,n}   \Vert \E ( Y_i | {\mathcal G}_{\ell-1}  )\Vert_p   \, .
\end{multline}
So, overall,  recalling that  $\eta_n = \sum_{i = 0}^n  \Vert \E (Y_i | {\mathcal G}_0 ) \Vert_p $, we get 
\begin{multline*}
U_{\ell,n} (p) \ll \max_{1 \leq \ell \leq n} | \alpha_{\ell ,n} |^{p-2}  \eta_n^{p-2}  \Big ( \sum_{i= \ell}^n  \alpha^2_{i,n}  (i-\ell +1)   \lambda_{i-\ell +1}  \\ +    \sum_{j= \ell}^{n }  \alpha^2_{n,j}   \sum_{u= [(j-\ell)/2]}^{n - j } \lambda_{u}  +   \eta_n \sum_{i=\ell}^n \alpha^2_{i,n}   \Vert \E ( Y_i | {\mathcal G}_{\ell-1}  )\Vert_p \Big ) \, .
\end{multline*}
Hence, setting 
\[
\Lambda_{i, \ell } = (i-\ell +1)   \lambda_{i-\ell +1} +   \sum_{u= [(i-\ell)/2]}^{n - i } \lambda_u + \eta_n \Vert \E ( Y_{i-\ell +1} | {\mathcal G}_{0}  )\Vert_p \, .
\]
we get 
\[
\sum_{ \ell =2}^n \frac{1}{ B^{(p-1)/2}_{ \ell , n} (a)  } U_{\ell,n} (p)  \ll  \max_{1 \leq \ell \leq n} | \alpha_{\ell ,n} |^{p-2}   \eta_n^{p-2}  \sum_{i= 1}^n   \frac{ \alpha^2_{i,n} }{ B^{(p-1)/2}_{ i , n} (a) }  \sum_{\ell =1}^i  \Lambda_{i, \ell }  \, .
\]
Since, for any  $i \leq n$, 
\[
\sum_{\ell =1}^i  \Lambda_{i, \ell }  \ll   \sum_{u=0}^n  ( (u+1) \lambda_u  +  \eta_n  \Vert \E ( Y_{u} | {\mathcal G}_{0}  )\Vert_p ) \leq  \Lambda_n  + \eta_n^2\, ,
\]
it follows that 
\beq \label{b1linUp}
\sum_{ \ell =2}^n \frac{1}{ B^{(p-1)/2}_{ \ell , n} (a)  } U_{\ell,n} (p)  \ll  \max_{1 \leq \ell \leq n} | \alpha_{\ell ,n} |^{p-2}   \eta_n^{p-2} (  \Lambda_n  + \eta_n^2)  \sum_{i= 1}^n   \frac{ \alpha^2_{i,n} }{ B^{(p-1)/2}_{ i , n} (a) }    \, .
\eeq
Let 
\[
a= \frac{   \max_{1 \leq k \leq n}  |\alpha_{k,n} |  \max(  \Vert Y_{0 } \Vert_2,  \sqrt{2 \pi  m})    + 2 \max_{1 \leq k \leq n-1} \Vert A_{k,n } \Vert_2   }{\max_{1 \leq k \leq n} \Vert d_{k,n}\Vert_2  }  \, ,
\]
where $m = \inf_{t \in [- \pi, \pi]} f_Y(t)$.  The decomposition \eqref{decmart1} entails that $a \geq 1$. On another hand, for any integer $\ell$ in  $ [1,n]$, 
\begin{align*}
B_{ \ell , n} (a)  & = \E (S_n-S_{\ell-1} - A_{\ell-1})^2 +a^2 \delta_n^2 = \E (S_n-S_{\ell-1})^2 - \E ( A_{\ell-1})^2 +a^2 \delta_n^2   \\
&  \geq   \Vert  S_n  -S_{\ell - 1} \Vert_2^2  +   \max_{1 \leq k \leq n}  |\alpha_{k,n} |^2  \max ( \Vert Y_{0 } \Vert^2_2, 2 \pi m   \big )   \, .
\end{align*}
But
\[
{\rm Var} (  S_n  -S_{\ell-1}) =  \int_{-\pi}^{\pi} \Big |  \sum_{k=\ell}^n  \alpha_{k ,n} {\rm e}^{{\rm i} t k} \Big |^2 f_Y(t) dt  \\
 \geq m \int_{-\pi}^{\pi} \Big |  \sum_{k=\ell}^n  \alpha_{k ,n} {\rm e}^{{\rm i} t k} \Big |^2 dt  = 2\pi m  \sum_{k=\ell}^n  \alpha_{k ,n}^2    \, .
\]
It follows that, for any integer $\ell$ in  $ [1,n]$, 
\beq  \label{upperboundBlalin*}
B_{ \ell , n} (a)  \geq 2 \pi m \Big (   \sum_{ i  = \ell }^n   \alpha_{i,n}^2  +  \max_{1 \leq k \leq n}  \alpha_{k,n}^2 \Big ) \, . 
\eeq
Starting from \eqref{b1linUp} and taking into account \eqref{upperboundBlalin*} and the fact that $m >0$, it follows that 
\begin{multline*}
\sum_{ \ell =2}^n \frac{1}{ B^{(p-1)/2}_{ \ell , n} (a)  } U_{\ell,n} (p)  \\ \ll  \max_{1 \leq \ell \leq n} | \alpha_{\ell ,n} |^{p-2}  \eta_n^{p-2}  (  \Lambda_n  + \eta_n^2)  \sum_{i= 1}^n   \frac{ \alpha^2_{i,n} }{ \Big (   \sum_{ j  = i }^n   \alpha_{j,n}^2  +  \max_{1 \leq k \leq n}  \alpha_{k,n}^2 \Big )^{(p-1)/2}}  \, .
\end{multline*}Hence proceeding as in the proof of Theorem \ref{Th1first}, we get 
\beq \label{b2linUp}
 \sum_{ \ell =2}^n \frac{1}{ B^{(p-1)/2}_{ \ell , n} (a)  } U_{\ell,n} (p)  \ll   \left\{
  \begin{aligned}
\max_{1 \leq \ell \leq n} | \alpha_{\ell ,n} |^{p-2}   \eta_n^{p-2} (  \Lambda_n  + \eta_n^2)   \Big (  \sum_{\ell= 1}^n \alpha^2_{\ell,n}  \Big  )^{(3-p)/2}  &  \ \text{ if $p\in (2,3) $}  \\
\max_{1 \leq \ell \leq n} | \alpha_{\ell ,n} |   \eta_n(  \Lambda_n  + \eta_n^2)  \log \Big ( m_n^{-1} \sum_{\ell= 1}^n \alpha^2_{\ell,n}  \Big  )   & \  \text{ if $p=3$.} \\
  \end{aligned}
\right. 
\eeq
On another hand, taking into account \eqref{b1dln} and proceeding as before we get
\beq \label{b3lindln}
 \sum_{ \ell =2}^n \frac{1}{ B^{(p-1)/2}_{ \ell , n} (a)  } \Vert d_{\ell ,n} \Vert^p_p \ll   \left\{
  \begin{aligned}
\max_{1 \leq \ell \leq n} | \alpha_{\ell ,n} |^{p-2}  \eta_n^{p}   \Big (  \sum_{\ell= 1}^n \alpha^2_{\ell,n}  \Big  )^{(3-p)/2}  & \ \text{ if $p\in (2,3) $}  \\
\max_{1 \leq \ell \leq n} | \alpha_{\ell ,n} |  \eta_n^{3}  \log \Big ( m_n^{-1} \sum_{\ell= 1}^n \alpha^2_{\ell,n}  \Big  )   &\  \text{ if $p=3 $.} \\
  \end{aligned}
\right. 
\eeq
Starting from \eqref{startlin} and taking into account   \eqref{boundterme4}, \eqref{b2linUp} and \eqref{b3lindln} together with the fact that
\[
a \delta_n \ll \max_{1 \leq \ell \leq n} | \alpha_{\ell ,n} |  \Big ( \sqrt{m} +  \sum_{k=0}^n \Vert \E (Y_k | {\mathcal G}_0) \Vert_2 \Big )  \, , 
\]
the corollary follows. 

\subsection{Proof of Corollary \ref{corapplilinear}}   The proof follows the lines of the proof of Corollary  \ref{corapplilinearcondonspectral}. The only difference is in the choice of $a$. We take here  
\[
a= \frac{   \max_{1 \leq k \leq n}  |\alpha_{k,n} |  \big ( \max ( \Vert Y_{0 } \Vert_2, \sqrt{2 \pi f_Y(0) } \big )   + 2 \max_{1 \leq k \leq n-1} \Vert A_{k,n } \Vert_2  + \sqrt{K(n)}}{\max_{1 \leq k \leq n} \Vert d_{k,n}\Vert_2  }  \, ,
\]
where $K(n) =  \big (  \sum_{k \geq 1}k^2|  \gamma_k |  \big )  \sum_{i=1}^{n+1} |\alpha_{i,n} - \alpha_{n,i-1}|^2$. 
Once again, the decomposition \eqref{decmart1} entails that $a \geq 1$. On another hand, 
\begin{align*}
B_{ \ell , n} (a)  & = \E (S_n-S_{\ell-1} - A_{\ell-1})^2 +a^2 \delta_n^2 = \E (S_n-S_{\ell-1})^2 - \E ( A_{\ell-1})^2 +a^2 \delta_n^2   \\
&  \geq   \Vert  S_n  -S_{\ell - 1} \Vert_2^2  +   \max_{1 \leq k \leq n}  |\alpha_{k,n} |^2 \big ( \max ( \Vert Y_{0 } \Vert^2_2, 2 \pi f_Y(0)  \big )  +  K(n) \, .
\end{align*}
But, setting ${\tilde \alpha}_u = \alpha_{u,n} $ if $u  \in [\ell , n]$ and $0$ otherwise, we get
\[
{\rm Var} (  S_n  -S_{\ell-1}) = \sum_{k \in {\mathbb Z}} \gamma_k \sum_{ i \in {\mathbb Z} } {\tilde \alpha}_i {\tilde \alpha}_{i+k} \\
= 2 \pi f_Y(0)  \sum_{ i  = \ell }^n   \alpha_{i,n}^2 - 2^{-1}  \sum_{k \in {\mathbb Z}} \gamma_k \sum_{ i \in {\mathbb Z} }  ( {\tilde \alpha}_i  - {\tilde \alpha}_{i+k} )^2  \, .
\]
Setting $K= \sum_{k \geq 1}k^2|  \gamma_k | $, it follows that 
\[
 \Vert  S_n  -S_{\ell-1} \Vert_2^2  +  K \sum_{i=1}^{n+1} |\alpha_{i,n} - \alpha_{i-1,n}|^2  \geq 2 \pi f_Y(0)  \sum_{ i  = \ell }^n   \alpha_{i,n}^2 \, , 
\]
implying that 
\[
B_{ \ell , n} (a)  \geq 2 \pi f_Y(0) \Big (   \sum_{ i  = \ell }^n   \alpha_{i,n}^2  +  \max_{1 \leq k \leq n}  \alpha_{k,n}^2 \Big ) \, . 
\]
Using the fact that  $f_Y(0) >0$, the rest of the proof is the same as that of Corollary \ref{corapplilinearcondonspectral}.

\subsection{Proof of Corollary \ref{corappli1}} We start as in the proof of Corollary \ref{corapplilinearcondonspectral} and use the notation introduced there. So we have the upper bound \eqref{startlin} with $p=3$. Recalling the notation $A_{k,n} = \E(S_n -  S_{k} | {\mathcal F}_k)$, we select 
\[
a= \frac{   \max_{1 \leq k \leq n} \Vert X_{k } \Vert_2  + 2 \max_{1 \leq k \leq n-1} \Vert A_{k,n } \Vert_2  }{\max_{1 \leq k \leq n} \Vert d_{k,n}\Vert_2  }  \, .
\]
The decomposition \eqref{decmart1} entails that  $a\geq 1$ and also that  
\[\sum_{i =k}^n \Vert d_{i,n} \Vert_2^2=  \Big \Vert \sum_{i =k}^n d_{i,n} \Big  \Vert_2^2 =  \big \Vert S_n - S_{k-1}   - A_{k-1,n}\big  \Vert_2^2 \, .
\]
It follows that \begin{align*}
B_{ k , n} (a)  & = \sum_{i =k}^n \Vert d_{i,n} \Vert_2^2 + a^2 \delta_n^2 \\& =  \Vert  S_n  -S_{k-1} \Vert_2^2 - \Vert A_{k-1,n}   \Vert_2^2 +  \big (  \max_{1 \leq k \leq n} \Vert X_{k } \Vert_2  + 2 \max_{1 \leq k \leq n-1} \Vert A_{k,n } \Vert_2  \big)^2  \\
& \geq   \Vert  S_n  -S_{k-1} \Vert_2^2  + \max_{1 \leq k \leq n} \Vert X_{k } \Vert_2^2 \, .
\end{align*}
Using $(H_2)$ and the fact that $C_n \geq 1$, we derive
\beq \label{majB}
\frac{1}{B_{ k , n} (a) } \leq \frac{C_n}{  \sum_{\ell=k}^n \Vert X_{\ell } \Vert_2^2 + \max_{1 \leq k \leq n} \Vert X_{k } \Vert_2^2} := \frac{C_n}{ {\widetilde V}_{k,n} } \, .
\eeq
On another hand, for any $1 \leq k \leq n-1$ and any $\eta >1/2$, by  the definition of the $\rho$-mixing coefficients, 
\begin{align*}
\Vert A_{k,n } \Vert^2_2  \leq \Big (  \sum_{\ell = k+1}^n \Vert   \E ( X_{\ell} | {\mathcal F}_k) \Vert_2 \Big )^2&  \ll  \sum_{\ell = k+1}^n (\ell -k )^{2 \eta} \Vert   \E ( X_{\ell} | {\mathcal F}_k) \Vert_2^2  \\
& \ll    \sum_{\ell = k+1}^n (\ell -k )^{2 \eta} \Vert  X_{\ell}  \Vert_2^2  \rho^2 (  \ell -k )  \, .
\end{align*}
According to $(H_1)$ we can take $\eta>1/2$ such that $ \sum_{\ell \geq 1} \ell^{2 \eta} \rho^2 (  \ell) < \infty$. Hence
\[
\Vert A_{k,n } \Vert_2 \ll  \max_{1 \leq k \leq n} \Vert X_k\Vert_2\, ,
\]
implying that 
\[
a\delta_n \ll \max_{1 \leq k \leq n} \Vert X_k\Vert_2 \, .
\]
On another hand, from decomposition \eqref{decmart1},
\[
\Vert d_{k,n} \Vert^3_3 \leq  9 \Big (  K_n \Vert X_{k} \Vert^2_2  +  \Vert A_{k,n } \Vert_3^3 + \Vert A_{k-1,n } \Vert_3^3 \Big ) \, .
\]
But, for any $1 \leq k \leq n-1$ and any $\eta >2/3$, by the definition of the $\rho$-mixing coefficients, 
\begin{multline*}
\Vert A_{k,n } \Vert^3_3  \leq \Big (  \sum_{\ell = k+1}^n \Vert   \E ( X_{\ell} | {\mathcal F}_k) \Vert_3 \Big )^3 \ll  \sum_{\ell = k+1}^n (\ell -k )^{3 \eta} \Vert   \E ( X_{\ell} | {\mathcal F}_k) \Vert_3^3  \\
\ll K_n  \sum_{\ell = k+1}^n (\ell -k )^{3 \eta} \Vert   \E ( X_{\ell} | {\mathcal F}_k) \Vert_2^2 \ll K_n  \sum_{\ell = k+1}^n (\ell -k )^{3 \eta} \Vert  X_{\ell}  \Vert_2^2  \rho^2 (  \ell -k )  \, .
\end{multline*}
So, overall,
\begin{align*}
 a \delta_n +      \sum_{k=1 }^n  \frac{ \E( |d_{k,n}|^3)}{ B_{ k+1 , n} (a)}  & \ll \max_{1 \leq k \leq n} \Vert X_k\Vert_2 + K_n C_n \sum_{k=1 }^n  \frac{\Vert X_k\Vert_2^2}{ {\widetilde V}_{k,n}}   \\ & \quad + K_nC_n
  \sum_{k=1 }^n  \frac{ \sum_{\ell = k}^n (\ell -k +1)^{3 \eta} \Vert  X_{\ell}  \Vert_2^2  \rho^2 (  \ell -k +1 ) }{ {\widetilde V}_{k,n}} 
  \\
 &   \ll K_n  + K_nC_n \sum_{k=1 }^n  \frac{\Vert X_k\Vert_2^2}{ {\widetilde V}_{k,n}}   \\ 
 & \quad + K_n C_n 
 \sum_{\ell = 1}^n   \frac{   \Vert  X_{\ell}  \Vert_2^2  }{ {\widetilde V}_{\ell,n}}  \sum_{k=1 }^\ell (\ell -k +1)^{3 \eta}  \rho^2 (  \ell -k +1 )
\end{align*}
According to $(H_1)$ we can take $\eta>2/3$ such that $ \sum_{\ell \geq 1} \ell^{3 \eta} \rho^2 (  \ell) < \infty$. Hence, it follows that 
\[
 a \delta_n +      \sum_{k=1 }^n  \frac{ \E( |d_{k,n}|^3)}{ B_{ k+1 , n} (a)} \ll K_n + K_n  C_n
 \sum_{\ell = 1}^n   \frac{   \Vert  X_{\ell}  \Vert_2^2  }{ {\widetilde V}_{\ell,n}} \, .
 \]
With similar arguments as those leading to \eqref{trivial2}, we get 
\[
 a \delta_n +      \sum_{k=1 }^n  \frac{ \E( |d_{k,n}|^3)}{ B_{ k+1 , n} (a)} \ll K_n + K_n C_n \log \Big  ( 1+ \sum_{k=1}^n \Vert X_k \Vert_2^2 \Big )   \, .
 \]

On another hand, we have
\[
U_{\ell,n} (3)  \leq 2   \Vert  d_{\ell -1,n}  \Vert_2  \Big \Vert \sum_{k = \ell}^n ( \E_{\ell-1} ( d_{k,n}^2) - \sigma_{k,n}^2)    \Big \Vert_2 \, . 
\]
To give an upper bound of this quantity we start from \eqref{p1martpart} with $p=4$. Note first that 
\begin{align*}
\Vert  \E^2_{\ell-1} (S_n-S_{\ell-1}) \Vert_2  & \leq 2 \sum_{i=\ell}^n \sum_{j=i}^n \Vert \E_{\ell-1} (X_i) \E_{\ell-1} (X_j) \Vert_2  \\
& \leq 2 \sum_{i=\ell}^n \sum_{j=i}^{2i-\ell} \Vert \E_{\ell-1} (X_i) X_j \Vert_2  +2 \sum_{i=\ell}^{n} \sum_{j=2i-\ell+1}^n \Vert X_i \E_{\ell-1} (X_j) \Vert_2  \, .
\end{align*}
Hence, by the definition of the $\rho$-mixing coefficients, we get 
\begin{equation} \label{p4martpart}
\Vert  \E^2_{\ell-1} (S_n-S_{\ell-1}) \Vert_2   \leq 4 K_n \sum_{i=\ell}^n  (i- \ell)  \Vert X_i \Vert_2 \rho (i - \ell)
\end{equation}
On another hand, by the definition of the $\rho$-mixing coefficients, we have: for $j \geq i \geq \ell$, 
\begin{equation} \label{p5martpart}
\Vert \E (X_iX_j | {\mathcal F}_{\ell-1}) -  \E (X_iX_j ) \Vert_2 \leq  \Vert X_iX_j \Vert_2 \rho (i- \ell) \leq  K_n  \Vert X_i \Vert_2 \rho (i- \ell)   \, , 
\end{equation}
and
\begin{equation} \label{p6martpart}
\Vert X_i  \E (X_j | {\mathcal F}_{i})  \Vert^2_2  = \E (  \E ( X_i^2X_j | {\mathcal F}_{i})  X_j)  \leq  K_n  \Vert X_i  \E (X_j | {\mathcal F}_{i})  \Vert_2  \Vert X_j \Vert_2 \rho (j- i ) \, .
\end{equation}
Hence starting from  \eqref{p1martpart} with $p=4$ and taking into account \eqref{p4martpart} and the upper bounds  \eqref{p2martpart} and \eqref{p3martpart} together with  \eqref{p5martpart} and  \eqref{p6martpart}, we derive
\[
U_{\ell,n} (3)  \ll  K_n  \Vert  d_{\ell -1,n}  \Vert_2 \sum_{i = \ell}^n  \Vert X_i \Vert_2 (i-\ell +1)  \rho ([i- \ell]/2 )  \, .
\]
Hence, taking into account $(H_1)$ and \eqref{majB}, 
\begin{align*}
\sum_{ \ell =2}^n \frac{1}{B_{ \ell , n} (a)  } U_{\ell,n} (3) & \ll K_n \sum_{ 2 \leq \ell  \leq i \leq n} \frac{ \Vert d_{\ell -1,n} \Vert_2  \Vert X_i \Vert_2}{B_{ \ell , n} (a)  } (i-\ell +1)  \rho ([i- \ell]/2 )  \\
& \ll  K_n \Big (  \sum_{ \ell =2}^n \frac{ \Vert d_{\ell -1,n} \Vert^2_2}{B_{ \ell , n} (a)  }  + \sum_{  i=2 }^n \frac{  \Vert X_i \Vert^2_2}{B_{ i , n} (a)  }   \Big )  \sum_{k =0}^n (k+1)  \rho (k/2)  \\
&  \ll  \Theta  K_n \Big (  \sum_{ \ell =2}^n \frac{ \Vert d_{\ell -1,n} \Vert^2_2}{B_{ \ell -1 , n} (a)  }  + C_n \sum_{  i=2 }^n \frac{  \Vert X_i \Vert^2_2}{{\widetilde V}_{ i, n} (a)  }   \Big ) \, ,
\end{align*}
since $B_{ \ell -1 , n} (a)  \leq 2 B_{ \ell  , n} (a) $. 
With similar arguments as those leading to \eqref{trivial2}, we get
\[
\sum_{ \ell =2}^n \frac{1}{B_{ \ell , n} (a)  } U_{\ell,n} (3) \ll  \Theta K_n  C_n \log  \Big ( 1+ \sum_{k=1}^n \Vert X_k \Vert_2^2 \Big ) \, .
\]
This ends the proof of the corollary since $ \sum_{k=1}^n \Vert X_k \Vert_2^2  \leq C_n V_n$.

\subsection{Proof of Corollary \ref{corsequential}}  As we shall see the result will use an approximation by a  ``reversed" martingale differences sequence. Hence, as a preliminary,  we first state the following fact: 

\begin{Fact}\textit{[Reversed martingale differences sequences] \label{commentreversed} Let $p \in (2,3]$.  Assume that $(d_n)_{n \in {\mathbb N}}$ is a real-valued  sequence of {\it  reversed } martingale differences  in ${\mathbb L}^p$ with respect to a  non-increasing sequence $({\mathcal G}_n)_{n \in {\mathbb N}}$ of $\sigma$-algebras. This means that for any integer $n$, $d_n$ is ${\mathcal G}_n$-adapted and ${\mathbb E} ( d_n | {\mathcal G}_{n+1} ) = 0$ a.s. Let $M_n = \sum_{k=1}^n d_k$. Note that $M_n = \sum_{k=1}^{n} \xi_{n,k}$ with $\xi_{n,k} = d_{n-k+1}$. Clearly $(\xi_{n,k})_{1 \leq k \leq n}$ is a sequence of martingale differences with respect to the increasing sequence $({\mathcal F}_{k,n})_k$ of $\sigma$-algebras with ${\mathcal F}_{k,n} = {\mathcal G}_{n-k+1}$. Hence, applying Proposition \ref{Mainprop}, it follows that \eqref{boundzetaprop} holds with ${\widetilde t}_{k, n} =\big (  \sum_{i=1}^{k -1} \E (d_i^2) + \delta^2 \big )^{1/2}$ replacing $t_{k,n}$, $d_k$ in place of $\xi_k$, ${\widetilde t}_{\ell +1, n}$ in place of $t_{\ell-1, n}$ and 
 \begin{equation} \label{deftildeU}
{ \tilde U}_{\ell,n} (p)  = \Big \Vert( |  d_{\ell+1} | \vee \sigma_{\ell+1} )^{p-2} \Big | \sum_{i= 1}^\ell ( \E ( d_i^2 | {\mathcal G}_{\ell+1}) - \sigma_i^2) \Big |   \Big \Vert_1
\end{equation}
in place of ${U}_{\ell,n} (p) $.  In particular, the following ``reversed" version of Theorem \ref{Th1first}  holds: setting ${\mathbb E}(d_i^2)=\sigma_i^2$ and $ \psi_n(t) = \sup_{1 \leq k \leq n }  \sigma_k^{-2} \E \inf ( t  \delta_n d_k^2 , |d_k|^3) $, 
there exist  positive constants $C_{r,p}$ depending on $( r, p)$ and $\kappa_r$ depending on $r$ such that for every positive integer $n$ and any  real $a\geq 1$,
\begin{multline} \label{boundzetarth1reverse}
\zeta_r(P_{M_{n}}, G_{V_n})  \leq C_{r,p}\Bigl (     \delta_n^r    \int_{a}^{\sqrt{ v_n (a) / \delta_n^2}} \frac{1}{x^{3-r}} dx   +   \delta_n^{r-1}\int_{a}^{\sqrt{ v_n(a) / \delta_n^2}} \frac{\psi_n ( \kappa_r  x ) }{x^{2-r}}dx  \\+ \sum_{ k =1}^{n-1} \frac{{ \tilde U}_{k,n} (p) }{ (a^2 \delta_n^2 + \sum_{i=1}^k \sigma_i^2)^{(p-r)/2}}  \Bigr ) + 4\sqrt{2}  a^r \delta_n^r
 \, .
\end{multline}}
\end{Fact}

We go back to the proof of Corollary \ref{corsequential}. Let ${\mathcal B}_n = \tau_n^{-1} \mathcal B$ and ${\tilde \phi}_k   =\phi_k - m(\phi_k (\tau_k)) $. As quoted by Conze and Raugi \cite{CR}, the following martingale-coboundary decomposition is valid: for any $n \in {\mathbb N}$, 
\beq \label{martdecdyn}
{\tilde \phi}_n = \psi_n -h_n +h_{n+1 } \circ T_{n+1} \, , 
\eeq
where $(d_n)_{n \geq 0}$ defined by $d_n =  \psi_n \circ \tau_n$ is a sequence of reversed martingale differences with respect to the filtration $({\mathcal B}_n)_{n \geq 0}$ and $(h_n)_{n \geq 0}$ is such that $m(h_n ( \tau_n))=0$, and there exists a positive constant $K$ such that $\sup_{n \geq 0} \Vert h_n \Vert_{\infty} \leq K$. 

Set $M_n = \sum_{k=1}^n d_n$ and $V(M_n) = \int_{X} M_n^2 (x) \, m(dx) =  \sum_{k=1}^n \int_{X} d_n^2 (x) \, m(dx) $.  We have
\[
W_1(P_{S_n}, G_{V_n})  \leq W_1(P_{S_n}, P_{M_n}) +W_1(P_{M_n}, G_{V(M_n)})  + W_1( G_{V(M_n)} ,G_{V_n}) \, .
\]
Using that $W_1( G_{V(M_n)} ,G_{V_n}))  \leq \big \vert  \sqrt{V(M_n)} - \sqrt{V_n} \big \vert  \leq \Vert  S_n - M_n\Vert_2$ and the martingale-coboundary decomposition  \eqref{martdecdyn}, it follows that 
\beq \label{dynW1-1}
W_1(P_{S_n}, G_{V_n})  \leq W_1(P_{M_n}, G_{V(M_n)}) + 4 \sup_{n \geq 0} \Vert h_n \Vert_{\infty}  \leq W_1(P_{M_n}, G_{V(M_n)})  + 4K\, .
\eeq
Since $\sup_{n \geq 0} \Vert d_n \Vert_{\infty} < \infty$, Corollary \ref{corsequential} will follow from Fact \ref{commentreversed} provided we can suitably  handle the quantities $ \Big \Vert \sum_{i= 1}^\ell ( \E ( d_i^2 | {\mathcal B}_{\ell+1}) - \E( d_i^2)  ) \Big \Vert_{1}$.  With this aim, note that  by \eqref{martdecdyn}, we have
\[
d_i^2   = {\tilde \phi}^2_i (\tau_i)  + 2 {\tilde \phi}_i (\tau_i)  ( h_i  (\tau_i) - h_{i+1 } ( \tau_{i+1} ) ) +  \big (  h_i  (\tau_i) - h_{i+1 } ( \tau_{i+1} ) \big )^2 \, ,
\]
implying that 
\begin{multline} \label{deccarre}
 \Vert \E ( d_i^2 | {\mathcal B}_{\ell+1}) - \E( d_i^2)   \Vert_{\infty} \leq   \Vert   \E (   {\tilde \phi}^2_i (\tau_i)  - m ( {\tilde \phi}^2_i (\tau_i) )   | {\mathcal B}_{\ell+1})   \Vert_{\infty} +   \Vert  \E (   h^2_i (\tau_i)  - m ( h^2_i (\tau_i) )   | {\mathcal B}_{\ell+1})   \Vert_{\infty} \\
+  \Vert  \E (   h^2_{i+1} (\tau_{i+1})  - m (  h^2_{i+1} (\tau_{i+1}) )   | {\mathcal B}_{\ell+1})  \Vert_{\infty}  + 2  \Vert  \E (   {\tilde \phi}_i (\tau_i)  h_i  (\tau_i)  - m (   {\tilde \phi}_i (\tau_i)  h_i  (\tau_i))  | {\mathcal B}_{\ell+1})    \Vert_{\infty} \\ + 2  \Vert   \E (  h_i (\tau_i)  h_{i+1} (\tau_{i+1})  - m (  h_i (\tau_i)  h_{i+1} (\tau_{i+1}) )   | {\mathcal B}_{\ell+1})   \Vert_{\infty}  \\
 +  2  \Vert  \E (   {\tilde \phi}_i (\tau_i)  h_{i+1} (\tau_{i+1})  - m (   {\tilde \phi}_i (\tau_i)   h_{i+1} (\tau_{i+1}) )  | {\mathcal B}_{\ell+1})  \Vert_{\infty}  \, .
\end{multline}
From Relations (1.8) and (1.10) in \cite{CR}, we get that for any function $f$  in ${\mathcal V}$ and  any $i \leq \ell $,  
\beq \label{espcondPerron}
  \E (  f (\tau_i)  - m ( f (\tau_i) )   | {\mathcal B}_{\ell+1})   =    \Big (   \frac{P_{\ell +1} \circ \cdots \circ P_{i+1} ({\tilde f}_i  \pi_{i}1 )  }{\pi_{\ell+1} 1} \Big )  \circ \tau_{\ell+1}  \, ,
\eeq
where ${\tilde f}_i = f - m ( f \pi_i 1 ) $. Hence taking into account  the properties  (DEC) and (MIN), we get that 
\begin{align*} 
 \Vert  \E (  f (\tau_i)  - m ( f (\tau_i) )   | {\mathcal B}_{\ell+1})   \Vert_{\infty} & \leq  \kappa  \delta^{-1} \Vert    P_{\ell +1} \circ \cdots \circ P_{i+1} ({\tilde f}_i  \pi_{i}1 )    \Vert_{v}  \\
&  \leq  \kappa  \delta^{-1}  C \gamma^{\ell+1 -i}   \Vert {\tilde f}_i  \pi_{i}1 \Vert_v \, .
\end{align*}
Hence, using Relation (3.10) in \cite{CR}, we get overall that there exists a positive constant $M$ such that, for any function $f$  in ${\mathcal V}$ and  any $i \leq \ell $,  
\beq \label{espcondPerronb1}
 \Vert  \E (  f (\tau_i)  - m ( f (\tau_i) )   | {\mathcal B}_{\ell+1})   \Vert_{\infty}  \leq M  \gamma^{\ell+1 -i}  \Vert f \Vert_v \, .
\eeq
Taking into account \eqref{espcondPerronb1}, it follows that the sum of the four first terms in the right-hand side of \eqref{deccarre} can be bounded by a positive constant times 
\beq \label{bounduniespcondPerron}
 \gamma^{\ell -i}  \Big (  \sup_{n \geq 0}  \Vert  h_n \Vert^2_v + \sup_{n \geq 0}   \Vert  \phi_n \Vert^2_v   \Big ) \, .
\eeq
To take care of the two last terms in \eqref{deccarre}, we shall use the following  fact:  for any functions $f$ and $g$ in ${\mathcal V}$, by using twice  \eqref{espcondPerron} and setting
\[
Q_{i+1}  f =  \frac{P_{i+1} ( f \pi_i 1 )  }{\pi_{i+1} 1} \, , 
\]
the following relation holds:   for any $i \leq \ell $,  
\begin{multline*}
\E (   f(\tau_i)  g (\tau_{i+1} )  | {\mathcal B}_{\ell+1})  = \E (  g (\tau_{i+1}  )  \E ( f(\tau_i)   | {\mathcal B}_{i+1})   | {\mathcal B}_{\ell+1})  \\
= \E \Big (g  \circ \tau_{i+1}    \Big (   \frac{P_{i+1} ( f \pi_i 1 )  }{\pi_{i+1} 1} \Big )  \circ \tau_{i+1}   \Big | {\mathcal B}_{\ell+1} \Big )  =   \Big (   \frac{P_{\ell +1} \circ \cdots \circ P_{i+2} ( g  Q_{i+1}  f  \pi_{i+1}1 )  }{\pi_{\ell+1} 1} \Big )  \circ \tau_{\ell+1}  \, .
\end{multline*}
Therefore, for any functions $f$ and $g$ in ${\mathcal V}$ and any $i \leq \ell $,
\begin{multline*}
\E (   f(\tau_i)  g (\tau_{i+1} )  -  m (  f(\tau_i)  g (\tau_{i+1} )  ) | {\mathcal B}_{\ell+1})  \\
 =  \Big (   \frac{P_{\ell +1} \circ \cdots \circ P_{i+2} (  ( g Q_{i+1}  f - m (g Q_{i+1}  f )  ) \pi_{i+1} 1 )    }{\pi_{\ell+1} 1} \Big )  \circ \tau_{\ell+1} \, .
\end{multline*}
Hence, taking into account the properties  (DEC) and (MIN), we get that for any $i \leq \ell $, 
\begin{align*}
 \Vert \E (   f(\tau_i)  g (\tau_{i+1} )  & -  m (  f(\tau_i)  g (\tau_{i+1} )  ) | {\mathcal B}_{\ell+1})  \Vert_{\infty}   \\
 & \leq  \kappa \delta^{-1}  \Vert P_{\ell +1} \circ \cdots \circ P_{i+2} (  ( g Q_{i+1}  f - m (g Q_{i+1}  f )  ) \pi_{i+1} 1 )   \Vert_v \\
&  \leq 
 \kappa  \delta^{-1}  C \gamma^{\ell -i} \Vert    ( g Q_{i+1}  f - m (g Q_{i+1}  f )  ) \pi_{i+1} 1 \Vert_v  \, .
\end{align*}
But
\begin{multline*}
\Vert    ( g Q_{i+1}  f - m (g Q_{i+1}  f )  ) \pi_{i+1} 1 \Vert_v \leq \Vert   (g Q_{i+1}  f )  \pi_{i+1} 1 \Vert_v+ \Vert  m (g Q_{i+1}  f )   \pi_{i+1} 1 \Vert_v  \\
\leq \Vert   g P_{i+1}  ( f  \pi_{i} 1 )  \Vert_v+ \Vert  g Q_{i+1}  f  \Vert_{\infty}  \Vert   \pi_{i+1} 1 \Vert_v  \leq  \kappa \Vert   g  \Vert_v \Vert P_{i+1}  ( f  \pi_{i} 1 ) \Vert_v+ \Vert  g Q_{i+1}  f  \Vert_{\infty}  \Vert   \pi_{i+1} 1 \Vert_v  \, .
\end{multline*}
By the property (DEC) we have $ \Vert P_{i+1} ( f \pi_i 1 ) \Vert_v \leq \kappa_3 \Vert f \Vert_v$ where $\kappa_3$ is a positive constant not depending on $i$ and on $f$.  On another hand, by the properties (DEC) and 
(MIN), we have 
\[
 \Vert  g Q_{i+1}  f  \Vert_{\infty}   \leq  \kappa \delta^{-1}   \Vert  g  \Vert_{\infty}  \Vert P_{i+1} ( f \pi_i 1 ) \Vert_v \leq \kappa_4   \Vert f \Vert_v \Vert g \Vert_v   \, ,
\]
where $\kappa_4$ is a positive constant not depending on $(i, f, g) $. 
So overall, there exists a positive constant $M$ such that, for any functions $f$  and $g$ in ${\mathcal V}$ and  any $i \leq \ell $,  
\beq \label{espcondPerronb2}
 \Vert  \E (  f (\tau_i)  g (\tau_{i+1} )  -  m (  f(\tau_i)  g (\tau_{i+1} )  )  | {\mathcal B}_{\ell+1})   \Vert_{\infty}  \leq M  \gamma^{\ell -i}  \Vert f \Vert_v \Vert g \Vert_v \, .
\eeq
Taking into account \eqref{espcondPerronb2}, it follows that the sum of the  two last terms in the right-hand side of \eqref{deccarre} can be bounded by a positive constant times the quantity \eqref{bounduniespcondPerron}.  So, overall, for any $i \leq \ell $,  
\beq \label{espcondcarre} 
\Big \Vert  |  d_{\ell+1}  |  \big |  \E ( d_i^2 | {\mathcal B}_{\ell+1}) - \E( d_i^2)  \big |   \Big \Vert_1  \ll \sup_{n \geq 0} \Vert d_n \Vert_{\infty} 
 \min ( \E( d_i^2)   ,  \gamma^{\ell -i}  )    \, .
\eeq
Therefore, recalling the notation \eqref{deftildeU} and setting $\delta_n^2= \max_{1 \leq i \leq n}  \E( d_i^2) $ and $a^2 = 1 + \delta_n^{-2}$, we get 
\[
\sum_{ \ell =1}^{n-1} \frac{{ \tilde U}_{\ell,n} (3) }{ a^2 \delta_n^2 + \sum_{k=1}^\ell \E( d_k^2) } \ll   \sum_{ \ell =1}^{n-1}  \sum_{i=1}^{\ell} \frac{  \min ( \E( d_i^2)   ,  \gamma^{\ell -i}  )  }{ 1 +  \delta_n^2 + \sum_{k=1}^i \E( d_k^2) } \, .
\]
Let $\alpha $ be a positive real and $\varphi_\alpha (\ell) = [\alpha \log (\ell)]$. Let $\ell_0 = \inf \{ \ell \geq 1 \, : \, \ell - \varphi_\alpha (\ell)  \geq 1 \}$. We then have
\begin{align*}
\sum_{ \ell =1}^{n-1} \frac{{ \tilde U}_{\ell,n} (3) }{ a^2 \delta_n^2 + \sum_{k=1}^\ell \E( d_k^2) } & \ll   \sum_{ \ell =1}^{n-1}  \sum_{i=1}^{\ell - \varphi_\alpha (\ell)}   \gamma^{\ell -i}   +  \sum_{ \ell =1}^{n-1}  \sum_{i = \ell - \varphi_\alpha (\ell) +1}^{\ell} \frac{  \E( d_i^2)     }{ 1 +  \delta_n^2 + \sum_{k=1}^i \E( d_k^2) } \\
& \ll   \sum_{ \ell =1}^{n-1}  ( 1 - \gamma)^{-1}  \gamma^{ \varphi_\alpha (\ell)}   +  (\log n ) \sum_{i =1}^{n-1} \frac{  \E( d_i^2)     }{ 1 +  \delta_n^2 + \sum_{k=1}^i \E( d_k^2) } 
 \, .
\end{align*}
Selecting $\alpha$ such that $\alpha \log ( 1/\gamma) >1$ and using similar arguments as those developed in Theorem \ref{Th1first}, it follows that
\[
\sum_{ \ell =1}^{n-1} \frac{{ \tilde U}_{\ell,n} (3) }{ a^2 \delta_n^2 + \sum_{k=1}^\ell \E( d_k^2) }  \ll  1  +  (\log n )   \log ( 1 + V(M_n)) \, . 
\]
Hence by taking into account this upper bound in \eqref{boundzetarth1reverse} (with $r=1$ and $p=3$), we derive that 
\begin{multline} \label{dynW1-2}
W_1(P_{M_n}, G_{V(M_n)}) + 4 \sup_{n \geq 0} \Vert h_n \Vert_{\infty} \\  \ll  1 + \sqrt{ \max_{1 \leq i \leq n}  \E( d_i^2)} +  \left(\max_{1 \leq i \leq n}  \Vert  d_i \Vert_{\infty} +  \log n \right )   \log ( 1 + V(M_n))  \, .
\end{multline}
Starting from \eqref{dynW1-1} and considering  \eqref{dynW1-2}  together with the fact that $\sup_{i \geq 1} \Vert d_i \Vert_{\infty} < \infty$ and that there exists a positive constant $B$ such that $V(M_n) \leq 2 V_n +B$, the result follows.

\end{document}